\documentclass[titlepage,draft,12pt]{article} 
\usepackage{amsfonts,latexsym,pstricks,pst-plot} 
\usepackage[a4paper]{geometry}

\newcommand{\re}{{\mathbb{R}}}

\newcommand{\ep}{\varepsilon}
\newcommand{\prf}{{\sc Proof.}$\;$}
\newcommand{\qed}{{\penalty 10000\mbox{$\quad\Box$}}}
\newcommand{\qone}{\mathcal{Q}_{1}}
\newcommand{\qtwo}{\mathcal{Q}_{2}}
\newcommand{\qthree}{\mathcal{Q}_{3}}
\newcommand{\qfour}{\mathcal{Q}_{4}}
\newcommand{\T}{\mathcal{T}}
\newcommand{\Tep}{\mathcal{T}_{\ep}}

\newcommand{\intp}{\mbox{\textup{Int}}_{P}}
\newcommand{\parbt}{\partial_{P}(\T)}
\newcommand{\uep}{u^{\ep}}
\newcommand{\fep}{f^{\ep}}

\newcommand{\D}{\mathcal{D}}
\newcommand{\kf}{\gamma_{2}}

\newtheorem{thm}{Theorem}[section]

\newtheorem{lemma}[thm]{Lemma}

\newtheorem{rmk}[thm]{Remark}

\title{An example of global classical solution for the Perona-Malik
equation} 

\author{Marina Ghisi\vspace{1ex}\\ {\normalsize Universit\`a degli
Studi di Pisa} \\{\normalsize Dipartimento di Matematica ``Leonida
Tonelli''}\\
{\normalsize PISA (Italy)}\\
{\normalsize e-mail: \texttt{ghisi@dm.unipi.it}}\and
Massimo Gobbino\vspace{1ex}\\ {\normalsize Universit\`a degli Studi di Pisa} 
\\{\normalsize Dipartimento di Matematica Applicata ``Ulisse Dini''}\\ 
{\normalsize 
 PISA (Italy)}\\  
{\normalsize e-mail: \texttt{m.gobbino@dma.unipi.it}}}

\date{}

\begin{document}
\maketitle
\begin{abstract}
	
	We consider the Cauchy problem for the Perona-Malik equation
	$$u_{t}=\mathrm{div}\left(\frac{\nabla u}{1+|\nabla
	u|^{2}}\right)$$
	in an open set $\Omega\subseteq\re^{n}$, with Neumann boundary
	conditions.
	
	It is well known that in the one-dimensional case this problem
	does not admit any global $C^{1}$ solution if the initial
	condition $u_{0}$ is transcritical, namely when $|\nabla
	u_{0}(x)|-1$ is a sign changing function in $\Omega$.
	
	In this paper we show that this result cannot be extended to
	higher dimension. We show indeed that for $n\geq 2$ the problem
	admits radial solutions of class $C^{2,1}$ with  a transcritical
	initial condition.
	
\vspace{1cm}

\noindent{\bf Mathematics Subject Classification 2000 (MSC2000):}
35K55, 35K65, 35B05.

\vspace{1cm} \noindent{\bf Key words:} Perona-Malik equation,
anisotropic diffusion, forward-backward parabolic equation, degenerate
parabolic equation, moving domains, overdetermined problem,
subsolutions and supersolutions.
\end{abstract}
 
\section{Introduction}

Let $\Omega\subseteq\re^{n}$ be  an open set. Let us consider the
Cauchy boundary value problem
\begin{eqnarray}
	u_{t}(x,t)-\mathrm{div}\left(\frac{\nabla u(x,t)}{1+|\nabla
	u(x,t)|^{2}}\right) & = & \makebox[3em][l]{0} \forall
	(x,t)\in\Omega\times[0,T), 
	\label{eq:PM-eq} \\
	\frac{\partial u}{\partial n}(x,t) & = & \makebox[3em][l]{0}
	\forall (x,t)\in\partial\Omega\times[0,T), 
	\label{eq:PM-nbc} \\
    u(x,0) & = & \makebox[3em][l]{$u_{0}(x)$}\forall x\in\Omega.
    \label{eq:PM-cauchy}
\end{eqnarray}

This problem was considered in \cite{PM} in the context of image
denoising.  In that framework $u_{0}$ is the grey level of a (noisy)
picture defined in a rectangle $\Omega\subseteq\re^{2}$, and $u(x,t)$
for a given $t>0$ should represent a restored image obtained by
smoothing the regions where $|\nabla u_{0}|<1$ (objects) and enhancing
the regions where $|\nabla u_{0}|>1$ (edges).

Equation (\ref{eq:PM-eq}) is the formal gradient flow of the
functional
$$PM(u):=\frac{1}{2}\int_{\Omega}^{}\log\left(1+|\nabla
u(x)|^{2}\right)dx.$$

The convex-concave behavior of the integrand makes (\ref{eq:PM-eq}) a
forward-backward partial differential equation of parabolic type, with
a forward (or subcritical) region where $|\nabla u(x,t)|<1$, and a
backward (or supercritical) region where $|\nabla u(x,t)|>1$.  If the
initial condition is subcritical, namely $|\nabla u_{0}(x)|<1$ for
every $x\in\Omega$, then the maximum principle guarantees that the
same is true for all subsequent times (see~\cite{kk}).  It follows
that in this case equation (\ref{eq:PM-eq}) is always forward
parabolic, and a smooth solution globally exists ($T=+\infty$).

Things are far more complicated when the initial condition is
\emph{transcritical}, namely when there are points $x\in\Omega$ where
$|\nabla u_{0}(x)|<1$, and points $x\in\Omega$ where $|\nabla
u_{0}(x)|>1$.  In this case the forward-backward character of the
equation makes the problem ill posed from the analytic point of view
(see \cite{kich}).  On the other hand, numerical computations exhibit
much more stability and show that the process has the desired
denoising effect on the initial condition.  This is usually referred
as the Perona-Malik paradox.  The mathematical understanding of this
phenomenology is still far away, despite of the considerable amount of
interest generated by this problem in the last fifteen years (see
\cite{BF,BNP,BNPT1,Z2,E1,E2,E3,guidotti,gl,Z1,Z3}).

In this paper we focus on classical solutions, namely solutions which
are at least of class $C^{1}$ or $C^{2,1}$.  Throughout this paper
$C^{2,1}$ denotes the standard parabolic space of functions with one
continuous derivative with respect to time, and two continuous
derivatives with respect to space variables.

As far as we know, in literature there are two main results for
classical solutions, and they are both limited to the one-dimensional
case.  The first known result was proved in \cite{GG-LocSol} and
concerns local-in-time solutions.  The result is that for $n=1$ the
set of initial data for which a local-in-time $C^{2,1}$ solution
exists is dense in $C^{1}(\Omega)$.  On the other hand, one cannot
expect existence of local classical solutions but for a dense set of
initial conditions (see \cite{kich}).

The second known result concerns global-in-time solutions, and it is a
nonexistence result.  It states that when $n=1$ there does not exist
any global-in-time solution of class $C^{1}$ with a transcritical
initial condition.  This was first proved in \cite{kk} with some
technical assumptions on $u_{0}$, and then in \cite{PMGlob} in full
generality.

Many attempts have been made to extend this result in higher
dimension.  The main point in the one-dimensional proof is the so
called \emph{persistence of supercritical regions}, namely the fact
that supercritical regions, if present for $t=0$, cannot disappear for
all subsequent times.  This important qualitative property of
solutions is also consistent with intuition, because it means that the
process does not destroy edges.

In dimension $n=1$ the persistence of supercritical regions follows
from the fact that the $L^{\infty}$ norm of $|\nabla u|$ in $\Omega$
is an increasing function of time.  Unfortunately this norm is
monotone for $n=1$, but in general not for $n\geq 2$, as shown by
Example~3 given in \cite[Section~4]{GG}.  A consequence of this
example is that the proof of the nonexistence result for $n=1$ cannot
be extended to $n\geq 2$.

In this paper we show that the nonexistence result itself cannot be
extended.  We show indeed that for $n\geq 2$ global-in-time classical
solutions \emph{do} exist, and also the persistence of supercritical
regions fails.  Our main result is the following (we state it for
$n=2$, but the generalization to any $n\geq 2$ is straightforward).

\begin{thm}\label{thm:main}
	Let $\Omega:=\{x\in\re^{2}:1\leq |x|\leq 5\}$ be an anulus.
	
	Then there exists $u\in C^{2,1}(\Omega\times[0,+\infty))$
	satisfying equation (\ref{eq:PM-eq}) with $T=+\infty$, the Neumann
	boundary condition (\ref{eq:PM-nbc}), and
	$$\left\{x\in\Omega:|\nabla u(x,0)|>1\right\}=
	\left\{x\in\Omega:2<|x|<4\right\}.$$
	
	Moreover there exists $t_{0}>0$ such that $|\nabla u(x,t)|<1$ for 
	every $(x,t)\in\Omega\times(t_{0},+\infty)$.
\end{thm}

By the way, before Theorem~\ref{thm:main} we didn't know any
nontrivial example of (even local-in-time) transcritical classical
solution in dimension $n\geq 2$.

In order to prove Theorem~\ref{thm:main} we can limit ourselves to
radial solutions.  Let $r:=|x|$ be a radial variable, and let us
consider radial solutions $u(r,t)$.  A simple computation (see
\cite{GG}) shows that equation (\ref{eq:PM-eq}) becomes the following
\begin{equation}
	u_{t}=\frac{1-u_{r}^{2}}{(1+u_{r}^{2})^{2}}u_{rr}+
	\frac{1}{r}\frac{u_{r}}{1+u_{r}^{2}},
	\label{eq:PM-radial}
\end{equation}
and the Neumann boundary condition (\ref{eq:PM-nbc}) in the anulus
becomes
\begin{equation}
	u_{r}(1,t)=0,
	\label{eq:nbc-l}
\end{equation}
\begin{equation}
	u_{r}(5,t)=0.
	\label{eq:nbc-r}
\end{equation}

Therefore Theorem~\ref{thm:main} is a consequence of the following
result.

\begin{thm}\label{thm:main-r}
	There exists $u\in C^{2,1}([1,5]\times[0,+\infty))$ satisfying
	equation (\ref{eq:PM-radial}) in the strip
	$[1,5]\times[0,+\infty)$, the Neumann boundary conditions
	(\ref{eq:nbc-l}) and (\ref{eq:nbc-r}) for every $t\geq 0$, and the
	estimates
	\begin{equation}
		0\leq u_{r}(r,0)<1
		\quad\quad\quad
		\forall r\in[1,2)\cup(4,5],
		\label{est:u0-sub}
	\end{equation}
	\begin{equation}
		u_{r}(r,0)>1
		\quad\quad\quad
		\forall r\in(2,4).
		\label{est:u0-super}
	\end{equation}
	
	Moreover there exists $t_{0}>0$ such that $|u_{r}(x,t)|<1$ for 
	every $(r,t)\in[1,5]\times(t_{0},+\infty)$.
\end{thm}

We recall that equation (\ref{eq:PM-radial}), without the second
summand in the right-hand side, is just the Perona-Malik equation in
dimension one.  A bureaucratic count of derivatives says that the
added term is a lower order term.  Nevertheless its influence on the
dynamic is enormous.  Due to that lower order term, the supercritical
region of our radial solution disappears in a finite time, in contrast
with the one dimensional case.  After the extinction of its supercritical
region, the solution becomes subcritical and has no more obstacles to
global existence.  For this reason what we need in the proof of
Theorem~\ref{thm:main-r} is to keep the solution alive and regular up
to this time.

We are afraid that the possible extinction of supercritical regions
makes the analytical study of the Perona-Malik equation even more
difficult in dimension $n\geq 2$.  We don't know whether this new
phenomenon had been observed before in numerical experiments.  We
leave to numerical analysts any discussion about its consequences on
the model and its practical applications.

We conclude with some comments on our main result.

Using an anulus instead of a ball should not be essential.
This choice spares us and the reader from the technicalities due to
the fact that equation (\ref{eq:PM-radial}) is singular for $r=0$.
The original equation however is not singular in the origin.
Therefore this singularity should be easily compensated by the first
Neumann boundary condition, which in a ball becomes $u_{r}(0,t)=0$.
For this reason it should not be difficult to find solutions with the 
same properties, but defined in a ball.

Concerning the nonlinearity, for simplicity we devoted this
introduction to the model case of the Perona-Malik equation.  On the
contrary, in the following sections we work in a more general setting.
What we actually do is to prove Theorem~\ref{thm:main} and
Theorem~\ref{thm:main-r} for the gradient flow of any integral
functional with nonconvex integrand.  We refer to the beginning of
section~\ref{sec:statements} for the details.

Finally, the reader could ask how special these solutions are.
Theorem~\ref{thm:main} indeed states the existence of just \emph{one}
such solution.  What we actually construct is a \emph{large class} of
such solutions.  In a word, they are not the fruit of some strange
pathology, but a common feature in dimension $n\geq 2$.  For example
in the proof of Theorem~\ref{thm:main-r}, after choosing $t_{0}$ small
enough, we have some freedom in the choice of the initial condition
$u(r,0)$ in the subcritical intervals $[1,2]$ and $[4,5]$, where we
only impose some inequalities and some compatibility conditions at the
endpoints (see section~\ref{sec:forward} for the details).  Our
construction then completes $u(r,0)$ in the supercritical interval
$[2,4]$ in such a way that the solution starting with that datum is
global and becomes subcritical for $t>t_{0}$.

This paper is organized as follows.  In section~\ref{sec:statements}
we reduce the construction of the required solution to the existence
of solutions in four suitable subdomains, in each one of which the
equation is either (degenerate) forward parabolic or (degenerate)
backward parabolic.  In section~\ref{sec:epsilon} we approximate the
degenerate problems with strictly parabolic problems depending on a
small parameter $\ep>0$, and we state several $\ep$-independent
estimates.  In section~\ref{sec:proofs} we prove these estimates.  In
section~\ref{sec:limit} we pass to the limit as $\ep\to 0^{+}$,
showing that the limits are solutions of the degenerate problems.
This completes the proof of Theorem~\ref{thm:main-r}, hence also of
Theorem~\ref{thm:main}.

\setcounter{equation}{0}
\section{The four subproblems}\label{sec:statements}

Let us introduce some notations.  Let $\D\subseteq\re^{2}$ be a
compact set.  The \emph{parabolic interior} of $\D$ is the set
$\intp(\D)$ of points $(r,t)\in\D$ for which there exists $\delta>0$
such that $[r-\delta,r+\delta]\times[t-\delta,t]\subseteq\D$. The
\emph{parabolic boundary} of $\D$ is the set
$\partial_{P}(\D):=\D\setminus\intp(\D)$.

Throughout this paper we assume that $\varphi\in C^{\infty}(\re)$ is
an even
function, hence in particular
\begin{equation}
	\varphi'(0)=\varphi'''(0)=0.
	\label{hp:phi'}
\end{equation}

Moreover we assume that
\begin{equation}
	\varphi''(\sigma)>0
	\hspace{2em}
	\forall\sigma\in[0,1),
	\label{hp:phi+}
\end{equation}
\begin{equation}
	\varphi''(1)=0,
	\label{hp:phi=}
\end{equation}
\begin{equation}
	\varphi''(\sigma)<0
	\hspace{2em}
	\forall\sigma\in(1,3],
	\label{hp:phi-}
\end{equation}
\begin{equation}
	\varphi'(3)>0.
	\label{hp:phi3}
\end{equation}

Of course the thresholds $\sigma=1$ and $\sigma=3$ can be replaced by
any pair of positive numbers $\sigma_{0}<\sigma_{1}$.

As a consequence of (\ref{hp:phi+}) and (\ref{hp:phi=}), or (\ref{hp:phi=})
and (\ref{hp:phi-}), we have that $\varphi'''(1)\leq 0$. As a
consequence of (\ref{hp:phi'}) through (\ref{hp:phi3}) we have also 
that $\varphi'(\sigma)>0$ for every $\sigma\in(0,3]$.

\begin{figure}[htbp]
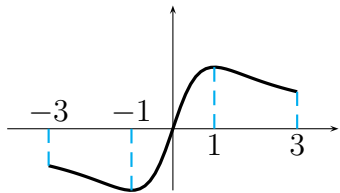

	\psset{unit=3ex}
	\centering
	\pspicture(-4,-1.5)(4,3)
	\psplot[linewidth=1.5\pslinewidth]{-3}{3}{x x 2 exp 1 add div 3 mul}
	\psline[linewidth=0.5\pslinewidth]{->}(-4,0)(4,0)
	\psline[linewidth=0.5\pslinewidth]{->}(0,-1.5)(0,3)
	\SpecialCoor
	\psline[linestyle=dashed, linecolor=cyan](1,0)(1,1.5)
	\psline[linestyle=dashed, linecolor=cyan](-1,0)(-1,-1.5)
	\psline[linestyle=dashed, linecolor=cyan](3,0)(3,0.9)
	\psline[linestyle=dashed, linecolor=cyan](-3,0)(-3,-0.9)
	\rput(1,-0.4){1}
	\rput(-1,0.4){$-1$}
	\rput(3,-0.4){3}
	\rput(-3,0.4){$-3$}
	\endpspicture

	\caption{typical graph of $\varphi'$}
	\label{fig:phi}
\end{figure}

Figure~\ref{fig:phi} shows the typical behavior of $\varphi'(\sigma)$.

We consider the following equation
\begin{equation}
	u_{t}(r,t)=\varphi''(u_{r}(r,t))u_{rr}(r,t)+
	\frac{\varphi'(u_{r}(r,t))}{r},
	\label{eq:PM-rad}
\end{equation}
which is the natural generalization of (\ref{eq:PM-radial}). 

It is easy to see that equation (\ref{eq:PM-rad}) reduces to
(\ref{eq:PM-radial}) when $\varphi(\sigma)=2^{-1}\log(1+\sigma^{2})$.
We believe and we hope that this generality simplifies the
presentation, and shows more clearly which properties of the
nonlinearity are essential in each step.

\subsection{The four regions}

In order to prove Theorem~\ref{thm:main-r} we divide the strip
$[1,5]\times [0,+\infty)$ into four regions.  To begin with, we fix
$t_{0}>0$, and we consider the functions
\begin{equation}
	\beta(t):=3-\sqrt{1-t/t_{0}},
	\hspace{3em}
	\gamma(t):=3+\sqrt{1-t/t_{0}},
	\label{defn:beta-gamma}
\end{equation}
defined for every $t\leq t_{0}$. Then we set
\begin{eqnarray}
	\qone & := & \left\{(r,t)\in\re^{2}:0\leq t\leq t_{0},\ 1\leq
	r\leq\beta(t)\right\},
	\label{defn:q1}  \\
	\qtwo & := & \left\{(r,t)\in\re^{2}:0\leq t\leq t_{0},\ \beta(t)\leq
	r\leq\gamma(t)\right\},
	\label{defn:q2}  \\
	\qthree & := & \left\{(r,t)\in\re^{2}:0\leq t\leq t_{0},\ \gamma(t)\leq
	r\leq 5\right\},
	\label{defn:q3}  \\
	\qfour & := & [1,5]\times[t_{0},+\infty).
	\label{defn:q4}
\end{eqnarray}

We also set
\begin{equation}
	\Gamma_{1}:=\left\{(\beta(t),t):t\in[0,t_{0}]\right\},
	\hspace{3em}
	\Gamma_{3}:=\left\{(\gamma(t),t):t\in[0,t_{0}]\right\}.
	\label{defn:gamma}
\end{equation}

\begin{figure}[tbp]
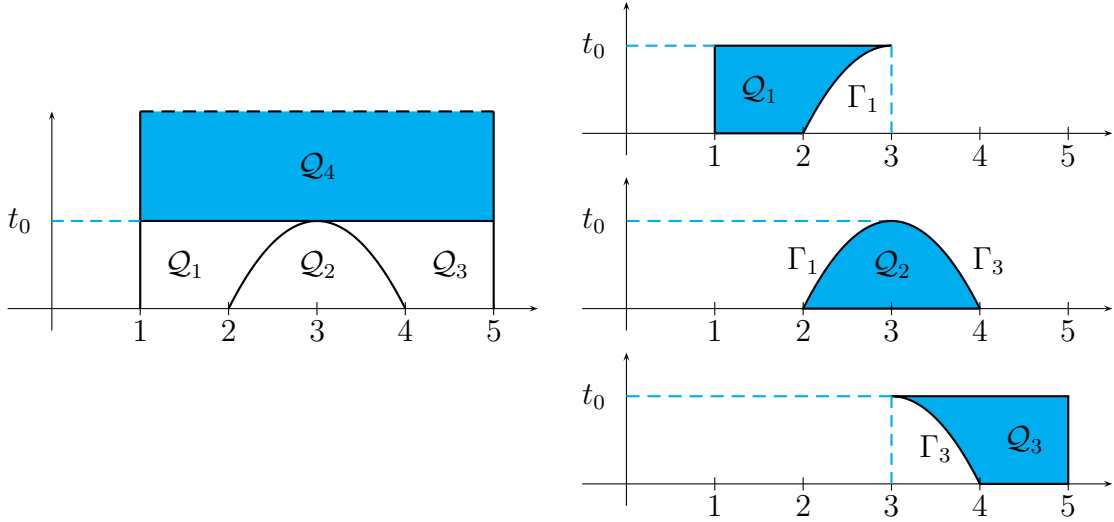

	\centering
\psset{yunit=3.2ex,xunit=6.4ex}
\pspicture(-1,-4.5)(12,7)
\psframe*[linecolor=cyan](1,2)(5,4.5)
\psplot{2}{4}{x 2 sub 4 x sub mul 2 mul}
\psline(1,2)(5,2)
\psline[linewidth=1\pslinewidth](1,0)(1,4.5)
\psline[linewidth=1\pslinewidth](5,0)(5,4.5)
\psline[linewidth=1\pslinewidth,linestyle=dashed](1,4.5)(5,4.5)
\psline[linewidth=0.5\pslinewidth]{->}(-0.5,0)(5.5,0)
\psline[linewidth=0.5\pslinewidth]{->}(0,-0.5)(0,4.5)
\psline[linewidth=0.5\pslinewidth](1,-0.15)(1,0.15)
\psline[linewidth=0.5\pslinewidth](2,-0.15)(2,0.15)
\psline[linewidth=0.5\pslinewidth](3,-0.15)(3,0.15)
\psline[linewidth=0.5\pslinewidth](4,-0.15)(4,0.15)
\psline[linewidth=0.5\pslinewidth](5,-0.15)(5,0.15)
\rput(1,-0.5){1}
\rput(2,-0.5){2}
\rput(3,-0.5){3}
\rput(4,-0.5){4}
\rput(5,-0.5){5}
\rput[l](-0.5,2){$t_{0}$}
\psline[linestyle=dashed, linecolor=cyan](0,2)(1,2)
\rput(3,1){$\qtwo$}
\rput(1.5,1){$\qone$}
\rput(4.5,1){$\qthree$}
\rput(3,3.25){$\qfour$}

\rput(6.5,4){
\pscustom[fillstyle=solid,fillcolor=cyan,linewidth=1\pslinewidth]{
\psline(3,2)(1,2)(1,0)(2,0)
\psplot{2}{3}{x 2 sub 4 x sub mul 2 mul}}
\psline[linewidth=0.5\pslinewidth]{->}(-0.5,0)(5.5,0)
\psline[linewidth=0.5\pslinewidth]{->}(0,-0.5)(0,3)
\psline[linewidth=0.5\pslinewidth](1,-0.15)(1,0.15)
\psline[linewidth=0.5\pslinewidth](2,-0.15)(2,0.15)
\psline[linewidth=0.5\pslinewidth](3,-0.15)(3,0.15)
\psline[linewidth=0.5\pslinewidth](4,-0.15)(4,0.15)
\psline[linewidth=0.5\pslinewidth](5,-0.15)(5,0.15)
\rput(1,-0.5){1}
\rput(2,-0.5){2}
\rput(3,-0.5){3}
\rput(4,-0.5){4}
\rput(5,-0.5){5}
\rput[l](-0.5,2){$t_{0}$}
\rput(1.5,1){$\qone$}
\rput[l](2.5,0.8){$\Gamma_{1}$}
\psline[linestyle=dashed, linecolor=cyan](0,2)(1,2)
\psline[linestyle=dashed, linecolor=cyan](3,0)(3,2)
}

\rput(6.5,0){
\psline[linestyle=dashed, linecolor=cyan](0,2)(3,2)
\pscustom[fillstyle=solid,fillcolor=cyan,linewidth=1\pslinewidth]{
\psline(4,0)(2,0)
\psplot{2}{4}{x 2 sub 4 x sub mul 2 mul}}
\psline[linewidth=0.5\pslinewidth]{->}(-0.5,0)(5.5,0)
\psline[linewidth=0.5\pslinewidth]{->}(0,-0.5)(0,3)
\psline[linewidth=0.5\pslinewidth](1,-0.15)(1,0.15)
\psline[linewidth=0.5\pslinewidth](2,-0.15)(2,0.15)
\psline[linewidth=0.5\pslinewidth](3,-0.15)(3,0.15)
\psline[linewidth=0.5\pslinewidth](4,-0.15)(4,0.15)
\psline[linewidth=0.5\pslinewidth](5,-0.15)(5,0.15)
\rput(1,-0.5){1}
\rput(2,-0.5){2}
\rput(3,-0.5){3}
\rput(4,-0.5){4}
\rput(5,-0.5){5}
\rput[l](-0.5,2){$t_{0}$}
\rput(3,1){$\qtwo$}
\rput(2,1.1){$\Gamma_{1}$}
\rput(4.1,1.1){$\Gamma_{3}$}
}

\rput(6.5,-4){
\pscustom[fillstyle=solid,fillcolor=cyan,linewidth=1\pslinewidth]{
\psline(4,0)(5,0)(5,2)(3,2)
\psplot{3}{4}{x 2 sub 4 x sub mul 2 mul}}
\psline[linewidth=0.5\pslinewidth]{->}(-0.5,0)(5.5,0)
\psline[linewidth=0.5\pslinewidth]{->}(0,-0.5)(0,3)
\psline[linewidth=0.5\pslinewidth](1,-0.15)(1,0.15)
\psline[linewidth=0.5\pslinewidth](2,-0.15)(2,0.15)
\psline[linewidth=0.5\pslinewidth](3,-0.15)(3,0.15)
\psline[linewidth=0.5\pslinewidth](4,-0.15)(4,0.15)
\psline[linewidth=0.5\pslinewidth](5,-0.15)(5,0.15)
\rput(1,-0.5){1}
\rput(2,-0.5){2}
\rput(3,-0.5){3}
\rput(4,-0.5){4}
\rput(5,-0.5){5}
\rput[l](-0.5,2){$t_{0}$}
\rput(4.5,1){$\qthree$}
\rput(3.5,0.8){$\Gamma_{3}$}
\psline[linestyle=dashed, linecolor=cyan](0,2)(3,2)
\psline[linestyle=dashed, linecolor=cyan](3,0)(3,2)
}
\endpspicture
	\caption{decomposition of $[1,5]\times[0,+\infty)$ into four
	subdomains}
	\label{fig:decomposition}
\end{figure}

These sets are represented in Figure~\ref{fig:decomposition}.
We finally consider the functions
\begin{equation}
	b(t):=\left\{
	\begin{array}{ll}
		\displaystyle{\min\left\{x\in\re:
		\varphi'''(1)x^{2}+\beta'(t)x-
		\frac{\varphi'(1)}{\beta^{2}(t)}=0\right\}} & 
		\mbox{if }t\in[0,t_{0}),\\
		0 & \mbox{if }t=t_{0},
	\end{array}
	\right.
	\label{defn:b}
\end{equation}
\begin{equation}
	c(t):=\left\{
	\begin{array}{ll}
		\displaystyle{\max\left\{x\in\re:
		\varphi'''(1)x^{2}+\gamma'(t)x-
		\frac{\varphi'(1)}{\gamma^{2}(t)}=0\right\}} & 
		\mbox{if }t\in[0,t_{0}),\\
		0 & \mbox{if }t=t_{0}.
	\end{array}
	\right.
	\label{defn:c}
\end{equation}

It is not difficult to see that $b(t)$ and $c(t)$ are well defined for
$t\in[0,t_{0}]$ provided that $t_{0}$ is small enough.  In
Lemma~\ref{lemma:b} and Lemma~\ref{lemma:c} below we prove some
properties of these functions, in particular their continuity.  In
Remark~\ref{rmk:heuristic} we explain the heuristic idea behind these
definitions.

In order to prove Theorem~\ref{thm:main-r} we need a solution of
(\ref{eq:PM-rad}) in the strip $[1,5]\times[0,+\infty)$, satisfying
the Neumann boundary conditions (\ref{eq:nbc-l}) and (\ref{eq:nbc-r})
for every $t\geq 0$, and estimates (\ref{est:u0-sub}) and
(\ref{est:u0-super}).  Our strategy is to construct this solution by
glueing together solutions of the same equation in the four
subdomains.  These solutions are of course required to fulfil the
Neumann boundary conditions (\ref{eq:nbc-l}) and (\ref{eq:nbc-r}).  In
order to be glued in a $C^{2,1}$ way along the common boundaries,
these solutions are also asked to satisfy the following conditions in
$\Gamma_{1}$
	\begin{equation}
		  \displaystyle{u(\beta(t),t)=u(3,t_{0})-3+\beta(t)-\varphi'(1)
		 \int_{t}^{t_{0}}\frac{1}{\beta(s)}\,ds}
		 \quad\quad\quad
		 \forall t\in[0,t_{0}], 
		\label{g1:u}  
	\end{equation}
	\begin{equation}
		  u_{r}(\beta(t),t)=1
		 \quad\quad\quad
		 \forall t\in[0,t_{0}],    
		\label{g1:ur}  
	\end{equation}
	\begin{equation}
		  u_{rr}(\beta(t),t)=b(t)
		 \quad\quad\quad
		 \forall t\in[0,t_{0}],    
		\label{g1:urr}  
	\end{equation}
and the following conditions in $\Gamma_{3}$
	\begin{equation}
		u(\gamma(t),t)=u(3,t_{0})-3+\gamma(t)-\varphi'(1)
		 \int_{t}^{t_{0}}\frac{1}{\gamma(s)}\,ds
		\quad\quad\quad
		\forall t\in[0,t_{0}],
		\label{g3:u}
	\end{equation}
	\begin{equation}
		u_{r}(\gamma(t),t)=1
		\quad\quad\quad
		 \forall t\in[0,t_{0}],
		\label{g3:ur}
	\end{equation}
	\begin{equation}
		u_{rr}(\gamma(t),t)=c(t)
		\quad\quad\quad
		\forall t\in[0,t_{0}].
		\label{g3:urr}
	\end{equation}

Moreover we arrange things in such a way that equation
(\ref{eq:PM-rad}) turns out to be (degenerate) backward parabolic in
$\qtwo$, and (degenerate) forward parabolic in the remaining three
regions.  Apparently all these conditions make the problem highly
overdetermined.  Nevertheless the degeneracy in $\Gamma_{1}$ and
$\Gamma_{3}$ allows us to fulfil all conditions.

\subsection{The four subproblems}

Let us state the existence results in the four regions.

\begin{thm}[Region $\qone$]\label{thm:q1}
	Let $\varphi\in C^{\infty}(\re)$ be a function satisfying
	(\ref{hp:phi'}), (\ref{hp:phi+}), and (\ref{hp:phi=}).
	Let us assume that $t_{0}>0$ is small enough.  Let $\beta(t)$,
	$\qone$, $\Gamma_{1}$, $b(t)$ be defined by (\ref{defn:beta-gamma}),
	(\ref{defn:q1}), (\ref{defn:gamma}), and (\ref{defn:b}),
	respectively.
	
	Then there exists $u\in C^{2,1}(\qone)$ satisfying equation
	(\ref{eq:PM-rad}) for every $(r,t)\in\qone$, the Neumann boundary
	condition (\ref{eq:nbc-l}) for every $t\in[0,t_{0}]$, and the
	boundary conditions (\ref{g1:u}) through (\ref{g1:urr}).
	
	Moreover $0\leq u_{r}(r,t)<1$ for every
	$(r,t)\in\qone\setminus\Gamma_{1}$.
\end{thm}

\begin{thm}[Region $\qtwo$]\label{thm:q2}
	Let $\varphi\in C^{\infty}(\re)$ be a function satisfying
	(\ref{hp:phi=}), (\ref{hp:phi-}), and (\ref{hp:phi3}).  Let us
	assume that $t_{0}>0$ is small enough.  Let $\beta(t)$,
	$\gamma(t)$, $\qtwo$, $\Gamma_{1}$, $\Gamma_{3}$, $b(t)$, $c(t)$
	be defined by (\ref{defn:beta-gamma}), (\ref{defn:q2}),
	(\ref{defn:gamma}), (\ref{defn:b}), and (\ref{defn:c}),
	respectively.
	
	Then there exists $u\in C^{2,1}(\qtwo)$ satisfying equation
	(\ref{eq:PM-rad}) for every $(r,t)\in\qtwo$ and the boundary
	conditions (\ref{g1:u}) through (\ref{g3:urr}).
	
	Moreover $u_{r}(r,t)>1$ for every
	$(r,t)\in\qtwo\setminus(\Gamma_{1}\cup\Gamma_{3}).$
\end{thm}

\begin{thm}[Region $\qthree$]\label{thm:q3}
	Let $\varphi\in C^{\infty}(\re)$ be a function satisfying
	(\ref{hp:phi'}), (\ref{hp:phi+}), and (\ref{hp:phi=}).
	Let us assume that $t_{0}>0$ is small enough.  Let $\gamma(t)$,
	$\qthree$, $\Gamma_{3}$, $c(t)$ be defined by
	(\ref{defn:beta-gamma}), (\ref{defn:q3}), (\ref{defn:gamma}), and
	(\ref{defn:c}), respectively.
	
	Then there exists $u\in C^{2,1}(\qthree)$ satisfying equation
	(\ref{eq:PM-rad}) for every $(r,t)\in\qthree$, the Neumann
	boundary condition (\ref{eq:nbc-r}) for every $t\in[0,t_{0}]$, and
	the boundary conditions (\ref{g3:u}) through (\ref{g3:urr}).
	
	Moreover $0\leq u_{r}(r,t)<1$ for every
	$(r,t)\in\qthree\setminus\Gamma_{3}$.
\end{thm}

\begin{thm}[Region $\qfour$]\label{thm:q4}
	Let $\varphi\in C^{\infty}(\re)$ be a function satisfying 
	(\ref{hp:phi'}), (\ref{hp:phi+}), and (\ref{hp:phi=}).
	Let $t_{0}>0$ be a real number, and let $\qfour$ be defined by
	(\ref{defn:q4}).  Let $u_{0}\in C^{2}([1,5])$ be a function such
	that
	\begin{equation}
		u_{0r}(1)=u_{0r}(5)=0,
		\label{hp:q4-n}
	\end{equation}
	\begin{equation}
		0\leq u_{0r}(r)<1
		\quad\quad\quad
		\forall r\in[1,3)\cup(3,5].
		\label{hp:q4-sp}
	\end{equation}
	
	Then there exists a unique function $u\in C^{2,1}(\qfour)$
	satisfying equation (\ref{eq:PM-rad}) for every $(r,t)\in\qfour$,
	the Neumann boundary condition (\ref{eq:nbc-l}) and
	(\ref{eq:nbc-r}) for every $t\geq t_{0}$, and the initial
	condition $u(r,t_{0})=u_{0}(r)$ for every $r\in[1,5]$.
	
	Moreover we have that $0\leq u_{r}(r,t)<1$ for every
	$(r,t)\in\qfour\setminus\{(3,t_{0})\}$.
\end{thm}

The idea of constructing a solution by glueing solutions in suitable
subdomains has already been used in \cite{GG-LocSol}. The main
difference is that in \cite{GG-LocSol} all subdomains are rectangles
(hence with \emph{fixed} endpoints), and the prescribed values of $u$,
$u_{r}$, $u_{rr}$ at the boundary do not depend on time.

Here the problem in $\qfour$ is a classical initial boundary value
problem in a fixed interval.  On the contrary, the problems in
$\qone$, $\qtwo$, $\qthree$ involve \emph{moving domains}, no initial
condition, but several \emph{time-dependent} boundary conditions in
the moving endpoints.

As for initial conditions, in $\qtwo$ the backward character of the
equation makes the solution completely determined by its values in
$\Gamma_{1}$ and $\Gamma_{3}$.  On the contrary, in $\qone$ and
$\qthree$ we have a lot of freedom in the choice of the initial
condition, and for this reason there are plenty of different
solutions.

Concerning the multiple boundary conditions, let us consider for
example the problem in $\qone$.  Once that an initial datum has been
chosen, the Neumann boundary conditions (\ref{eq:nbc-l}) and
(\ref{g1:ur}) are enough to determine uniquely the solution.  The
surprising aspect is that this solution satisfies also (\ref{g1:u})
and (\ref{g1:urr}), independently on the initial condition!  In the
following two remarks we show how the degeneracy of the equation in
$\Gamma_{1}$ makes this possible.

\begin{rmk}\label{rmk:heuristic}
	\begin{em}
		Let $u$ be a solution of equation (\ref{eq:PM-rad}) in $\qone$
		or $\qtwo$, with Neumann boundary condition (\ref{g1:ur}).  If
		$u$ is of class $C^{2,1}$ we can compute
		\begin{eqnarray*}
			\frac{\mbox{d}}{\mbox{d}t}\left[u(\beta(t),t)\right] & = & 
			\beta'(t)u_{r}(\beta(t),t)+u_{t}(\beta(t),t)  \\
			 & = & \beta'(t)u_{r}(\beta(t),t)+
			 \varphi''(u_{r}(\beta(t),t))u_{rr}(\beta(t),t)+
			 \frac{\varphi'(u_{r}(\beta(t),t))}{\beta(t)}.
		\end{eqnarray*}
		
		From condition (\ref{g1:ur}) we have
		therefore that
		$$\frac{\mbox{d}}{\mbox{d}t}\left[u(\beta(t),t)\right]=
		\beta'(t)+\frac{\varphi'(1)}{\beta(t)}.$$
		
		Integrating this equality in $[t,t_{0}]$ we obtain
		(\ref{g1:u}).  This is actually a proof of (\ref{g1:u}).  The
		proof of (\ref{g3:u}) is analogous.
	\end{em}
\end{rmk}

\begin{rmk}
	\begin{em}
		Let $u$ be a solution of equation (\ref{eq:PM-rad}) in $\qone$
		or $\qtwo$, with Neumann boundary condition (\ref{g1:ur}).  Let
		us assume that $u$ is smooth enough, and let us compute the
		time derivative of (\ref{g1:ur}).  We obtain that
		\begin{eqnarray*}
			0\ =\
			\frac{\mbox{d}}{\mbox{d}t}\left[u_{r}(\beta(t),t)\right] &
			= & \beta'(t)u_{rr}(\beta(t),t)+u_{rt}(\beta(t),t) \\
			 & = & \beta'(t)u_{rr}+ \varphi''(u_{r})u_{rrr}+
			\varphi'''(u_{r})u_{rr}^{2}+
			\frac{\varphi''(u_{r})u_{rr}}{\beta(t)}-
			\frac{\varphi'(u_{r})}{\beta^{2}(t)},
		\end{eqnarray*}
		where in the last line all the derivatives of $u$ are computed
		in $(\beta(t),t)$.  Exploiting the Neumann boundary condition
		(\ref{g1:ur}) we have therefore that
		\begin{equation}
			\beta'(t)u_{rr}(\beta(t),t)+
			\varphi'''(1)u_{rr}^{2}(\beta(t),t)-
			\frac{\varphi'(1)}{\beta^{2}(t)}=0.
			\label{eq:urr}
		\end{equation}
		
		This means that $u_{rr}(\beta(t),t)$ is a solution of the
		equation defining $b(t)$ in (\ref{defn:b}). If for $t=0$
		we have that $u_{rr}(\beta(0),0)$ is the smallest solution of 
		the equation, namely $b(0)$, then for all subsequent times 
		$u_{rr}(\beta(t),t)$ remains the smallest solution, namely
		$b(t)$. The choice of the smallest solution is due to the fact
		that it tends to $0$ as $t\to t_{0}^{-}$, while the other
		solution diverges to $+\infty$.
		
		We point out that this simple argument is a heuristic
		justification of (\ref{g1:urr}), but it is by no means a
		proof.  In deriving (\ref{eq:urr}) we used indeed that terms
		such as $\varphi''(u_{r})u_{rrr}$ vanish at the moving
		boundary.  This requires some assumption on $u_{rrr}$, which
		is beyond the $C^{2,1}$ regularity.
	\end{em}
\end{rmk}

\subsection{Proof of Theorem~\ref{thm:main-r}}

Let us take $t_{0}>0$ small enough in order to apply
Theorem~\ref{thm:q1}, Theorem~\ref{thm:q2}, and Theorem~\ref{thm:q3}.
Thus we obtain a solution of equation (\ref{eq:PM-rad}) in each one of
the regions $\qone$, $\qtwo$, $\qthree$.  All these solutions are
defined up to additive constants.  We can therefore assume, without
loss of generality, that they coincide in the common point
$(3,t_{0})$.

We claim that the three solutions glue together in order to give a
solution of class $C^{2,1}$ in the rectangle $[1,5]\times[0,t_{0}]$.

Let us examine indeed the solutions in $\qone$ and $\qtwo$.  If they
coincide in $(3,t_{0})$, then they coincide in the whole $\Gamma_{1}$
because they both satisfy (\ref{g1:u}).  Also their first and second
space derivatives coincide in $\Gamma_{1}$ because they both satisfy
(\ref{g1:ur}) and (\ref{g1:urr}).  So these two solutions glue in a
$C^{2,1}$ way.

The same is true for the solutions in $\qtwo$ and $\qthree$.  Note
that in the common point $(3,t_{0})$ we have that
$u_{rr}(3,t_{0})=b(t_{0})=c(t_{0})=0$.

It remains to extend $u$ to $[1,5]\times[t_{0},+\infty)$.  To this
end, we apply Theorem~\ref{thm:q4} using as ``initial'' condition the
trace at $t=t_{0}$ of the solution we have just glued in
$[1,5]\times[0,t_{0}]$.  From Theorem~\ref{thm:q1} and
Theorem~\ref{thm:q3} it is clear that this trace satisfies assumptions
(\ref{hp:q4-n}) and (\ref{hp:q4-sp}).  This completes the proof of
Theorem~\ref{thm:main-r}.  \qed

\setcounter{equation}{0}
\section{Approximating problems}\label{sec:epsilon}

In the previous section the proof of our main result has been reduced
to the proof of Theorem~\ref{thm:q1} through Theorem~\ref{thm:q4}.  In
this section we present our approach to Theorem~\ref{thm:q1} and
Theorem~\ref{thm:q2}.  We skip Theorem~\ref{thm:q3} because it is
symmetric to Theorem~\ref{thm:q1}, and we skip Theorem~\ref{thm:q4}
because it concerns a quite standard (non overdetermined) initial
boundary value problem.

The first thing to do is to choose $t_{0}$.  Let us begin by
considering the following constants depending only on $\varphi$
\begin{eqnarray}
	 & \gamma_{0}:=3\varphi'(1)+5, & 
	\label{defn:gamma0}  \\
	 & \gamma_{1}:=5\varphi'(1)+100, & 
	\label{defn:gamma1}  \\
	 & \kf:=\max\left\{
	|\varphi'(\sigma)|+|\varphi''(\sigma)|+
	|\varphi'''(\sigma)|+|\varphi^{IV}(\sigma)|:
	\sigma\in[0,3]\right\}. & 
	\label{defn:gamma2}
\end{eqnarray}

Let us choose $t_{0}\in(0,1)$ satisfying the following inequalities
\begin{equation}
 	t_{0}\leq\frac{1}{4\left[\varphi'(1)\right]^{2}},
 	\quad\quad\quad
	t_{0}\leq\frac{3}{2500\kf},
	\quad\quad\quad
	t_{0}\leq\frac{1}{96(\gamma_{1}+1)^{4}\kf},
	\label{hp:t01}
\end{equation}
\begin{equation}
	t_{0}\leq\frac{1}{(20\gamma_{0}^{2}+28\gamma_{0}+9)\kf},
	\quad\quad\quad
	t_{0}\leq\frac{1}{(12\gamma_{0}+14)\kf}.
	\label{hp:t02}
\end{equation}

We stated these conditions on $t_{0}$ as they are required throughout
the proofs. It is not difficult to see that the third inequality in
(\ref{hp:t01}) implies the remaining four.

\subsection{The forward problem in a moving domain}\label{sec:forward}

Let us consider the problem in $\qone$.  Let us choose $u_{0}\in
C^{\infty}([1,2])$ satisfying the compatibility conditions
\begin{equation}
	u_{0r}(1)=0,
	\quad\quad
	u_{0r}(2)=1,
	\quad\quad
	u_{0rr}(2)=b(0),
	\label{hp:u0-comp}
\end{equation}
and the following inequalities
\begin{eqnarray}
	 & 0\leq u_{0r}(r)<1
	\hspace{3em}
	\forall r\in(1,2), & 
	\label{hp:u0r}  \\
	 & |u_{0rr}(r)|<10
	\hspace{3em}
	\forall r\in[1,2], & 
	\label{hp:u0rr}  \\
	 & |u_{0rrr}(r)|<10
	\hspace{3em}
	\forall r\in[1,2]. & 
	\label{hp:u0rrr}
\end{eqnarray}

We remind that $b(0)$ depends on $t_{0}$, and for this reason the
choice of $u_{0}$ depends on the choice of $t_{0}$.

Now we approximate the problem in $\qone$ with a strictly parabolic
problem depending on a small parameter $\ep$.

\begin{thm}\label{thm:Qep}
	Let $\varphi\in C^{\infty}(\re)$ be a function satisfying
	(\ref{hp:phi'}), (\ref{hp:phi+}), (\ref{hp:phi=}).  Let us assume
	that $t_{0}>0$ satisfies (\ref{hp:t01}) and (\ref{hp:t02}).  Let
	$\beta(t)$, $\qone$, $b(t)$, be defined by
	(\ref{defn:beta-gamma}), (\ref{defn:q1}), and (\ref{defn:b}),
	respectively.  Let $u_{0}\in C^{\infty}([1,2])$ be a function
	satisfying (\ref{hp:u0-comp}) through (\ref{hp:u0rrr}), and let
	$\ep\in(0,1)$.
	
    Then there exists a unique function $\uep\in C^{2,1}(\qone)$
	satisfying equation (\ref{eq:PM-rad}) for every $(r,t)\in\qone$,
	the Neumann boundary condition (\ref{eq:nbc-l}) for every
	$t\in[0,t_{0}]$, the further Neumann boundary condition
	\begin{equation}
		\uep_{r}(\beta(t),t)=1-\ep
		\quad\quad\quad
		 \forall t\in[0,t_{0}],
		\label{th:Q-nr}
	\end{equation}
	and the initial condition
	$$\uep(r,0)=(1-\ep)u_{0}(r)
		\quad\quad\quad
		 \forall r\in[1,2].$$

    Moreover $\uep$ satisfies the following estimates independent on
	$\ep$.
    \begin{enumerate}
		\item \emph{Maximum principle for space derivatives}.  We have
		that
		\begin{equation}
			0\leq\uep_{r}(r,t)\leq 1-\ep
			\quad\quad\quad
			\forall (r,t)\in\qone.
			\label{est:Q-ur}
		\end{equation}
		
		\item \emph{Uniform strict parabolicity in the interior}.  For
		every $\delta\in(0,1)$ there exists a constant
		$M_{1}=M_{1}(\delta)$ such
		that
		\begin{equation}
			\uep_{r}(r,t)\leq M_{1}<1
			\hspace{3em}
			\forall t\in[0,t_{0}],\ \forall r\in[1,\beta(t)-\delta].
			\label{est:Q-unif-farab}
		\end{equation}
	
		As a consequence, there exists a constant
		$M_{2}=M_{2}(\delta)>0$ such that
		\begin{equation}
			\varphi''(\uep_{r}(r,t))\geq M_{2} 
			\hspace{3em}
			\forall t\in[0,t_{0}],\ \forall r\in[1,\beta(t)-\delta].
			\label{est:Q-unif-parab}
		\end{equation}
	
		\item \emph{Estimate on second derivatives at the fixed
		boundary}.  We have that
		\begin{equation}
	    	0\leq\uep_{rr}(1,t)\leq 100
			\hspace{3em}
	        \forall t\in[0,t_{0}].
	    	\label{est:Q-ueprr-l}
		\end{equation}

		\item \emph{Estimate on second derivatives at the moving
		boundary}.  We have that
		\begin{equation}
	    	\left|\uep_{rr}(\beta(t),t)-b(t)\right|
			\leq \ep
			\hspace{3em}
	        \forall t\in[0,t_{0}].
	    	\label{est:Q-ueprr}
		\end{equation}

		\item \emph{Global estimate on second derivatives}.  There
		exists a constant $M_{3}$ such that
		\begin{equation}
	    	\left|\uep_{rr}(r,t)-b(t)\right|
			\leq \ep+M_{3}\left|r-\beta(t)\right|
			\hspace{3em}
	        \forall (r,t)\in\qone.
	    	\label{est:Q-ueprr-n}
		\end{equation}
		
		As a consequence, there exist constants $M_{4}$ and $M_{5}$
		such that
		\begin{equation}
			|\uep_{rr}(r,t)|\leq M_{4}
			\quad\quad\quad
			\forall (r,t)\in\qone,
			\label{est:ueprr-glob}
		\end{equation}
		\begin{equation}
			|\uep_{t}(r,t)|\leq M_{5}
			\quad\quad\quad
			\forall (r,t)\in\qone.
			\label{est:uept-glob}
		\end{equation}
    
		\item \emph{Integral estimates}.  There exists a
		constant $M_{6}$ such that
		\begin{equation}
			\int_{\qone}\varphi''(\uep_{r}(r,t))
			[\uep_{rt}(r,t)]^{2}\,dr\,dt\leq M_{6}.
			\label{est:ueprt-int}
		\end{equation}
		
		In addition, for every $\delta\in(0,1)$ there exists a
		constant $M_{7}=M_{7}(\delta)$ such that
		\begin{equation}
			\int_{1}^{\beta(t)-\delta}[\uep_{rt}(r,t)]^{2}\,dr  \leq 
	    	M_{7} 
			\hspace{3em} 
			\forall t\in[0,t_{0}],
			\label{est:ueprt}
		\end{equation}
		\begin{equation}
			\int_{1}^{\beta(t)-\delta}[\uep_{rrr}(r,t)]^{2}\,dr \leq
			M_{7}
			\hspace{3em} 
			\forall t\in[0,t_{0}],
			\label{est:ueprrr-int}
		\end{equation}
		\begin{equation}
			\int_{0}^{t_{0}}\int_{1}^{\beta(t)-\delta}
	    	[\uep_{rrt}(r,t)]^{2}\,dr\,dt  \leq  M_{7}.
	    	\label{est:ueprrt-int}
		\end{equation}
\end{enumerate}
\end{thm}

\subsection{The backward problem in a moving domain}\label{sec:back}

The backward problem in $\qtwo$ can be transformed in a
forward parabolic problem by reversing the time.  To this end we
consider the domain
\begin{equation}
	\T:=\left\{(r,t)\in\re^{2}:0\leq t\leq t_{0},\ \beta(t_{0}-t)\leq
	r\leq\gamma(t_{0}-t)\right\}.
	\label{defn:T}
\end{equation}

The domain $\T$ is just $\qtwo$ upside-down, while its parabolic
boundary $\partial_{P}\T$ is
just $\Gamma_{1}\cup\Gamma_{3}$ upside-down (see Figure~\ref{fig:T}).
\begin{figure}[htbp]
	\centering
\psset{unit=3.8ex}
\pspicture(-1,-1)(6,2.5)
\psline[linestyle=dashed, linecolor=cyan](0,2)(3,2)
\pscustom[fillstyle=solid,fillcolor=green,linewidth=1\pslinewidth]{
\psline(4,0)(2,0)
\psplot{2}{4}{x 2 sub 4 x sub mul 2 mul}}
\psline[linewidth=0.5\pslinewidth]{->}(-0.5,0)(6,0)
\psline[linewidth=0.5\pslinewidth]{->}(0,-0.5)(0,3)
\psline[linewidth=0.5\pslinewidth](1,-0.15)(1,0.15)
\psline[linewidth=0.5\pslinewidth](2,-0.15)(2,0.15)
\psline[linewidth=0.5\pslinewidth](3,-0.15)(3,0.15)
\psline[linewidth=0.5\pslinewidth](4,-0.15)(4,0.15)
\psline[linewidth=0.5\pslinewidth](5,-0.15)(5,0.15)
\rput(1,-0.5){1}
\rput(2,-0.5){2}
\rput(3,-0.5){3}
\rput(4,-0.5){4}
\rput(5,-0.5){5}
\rput(-0.5,2){$t_{0}$}
\rput(3,0.8){$\qtwo$}
\endpspicture
\hfill
\pspicture(-1,-1)(6,3)
\psline[linestyle=dashed, linecolor=cyan](0,2)(2,2)
\pscustom[fillstyle=solid,fillcolor=green,linewidth=1\pslinewidth]{
\psline(4,2)(2,2)
\psplot{2}{4}{2 x 2 sub 4 x sub mul 2 mul sub}}
\psline[linewidth=0.5\pslinewidth]{->}(-0.5,0)(6,0)
\psline[linewidth=0.5\pslinewidth]{->}(0,-0.5)(0,3)
\psline[linewidth=0.5\pslinewidth](1,-0.15)(1,0.15)
\psline[linewidth=0.5\pslinewidth](2,-0.15)(2,0.15)
\psline[linewidth=0.5\pslinewidth](3,-0.15)(3,0.15)
\psline[linewidth=0.5\pslinewidth](4,-0.15)(4,0.15)
\psline[linewidth=0.5\pslinewidth](5,-0.15)(5,0.15)
\rput(1,-0.5){1}
\rput(2,-0.5){2}
\rput(3,-0.5){3}
\rput(4,-0.5){4}
\rput(5,-0.5){5}
\rput(-0.5,2){$t_{0}$}
\rput(3,1.1){$\T$}
\endpspicture
\hfill
\pspicture(-1,-1)(6,3)
\psline[linestyle=dashed, linecolor=cyan](0,2)(2,2)
\psline[linestyle=dashed, linecolor=cyan](0,0.32)(2.6,0.32)
\psplot[linewidth=0.5\pslinewidth]{2}{4}{2 x 2 sub 4 x sub mul 2 mul sub}
\pscustom[fillstyle=solid,fillcolor=green,linewidth=1\pslinewidth]{
\psline(4,2)(2,2)
\psplot{2}{2.6}{2 x 2 sub 4 x sub mul 2 mul sub}
\psline(2.6,0.32)(3.4,0.32)
\psplot{3.4}{4}{2 x 2 sub 4 x sub mul 2 mul sub}}
\psline[linewidth=0.5\pslinewidth]{->}(-0.5,0)(6,0)
\psline[linewidth=0.5\pslinewidth]{->}(0,-0.5)(0,3)
\psline[linewidth=0.5\pslinewidth](1,-0.15)(1,0.15)
\psline[linewidth=0.5\pslinewidth](2,-0.15)(2,0.15)
\psline[linewidth=0.5\pslinewidth](3,-0.15)(3,0.15)
\psline[linewidth=0.5\pslinewidth](4,-0.15)(4,0.15)
\psline[linewidth=0.5\pslinewidth](5,-0.15)(5,0.15)
\rput(1,-0.5){1}
\rput(2,-0.5){2}
\rput(3,-0.5){3}
\rput(4,-0.5){4}
\rput(5,-0.5){5}
\rput(-0.5,2){$t_{0}$}
\rput(3,1.2){$\Tep$}
\rput(-0.4,0.32){$\varepsilon$}
\endpspicture
	
	\caption{domains $\qtwo$, $\T$, and $\Tep$}
	\label{fig:T}
\end{figure}

With this notation proving Theorem~\ref{thm:q2} is equivalent to
proving the following result (note in particular the minus signs in
the right-hand side of (\ref{eq:T-eq})).

\begin{thm}[Region $\T$]\label{thm:T}
	Let $\varphi\in C^{\infty}(\re)$ be a function satisfying
	(\ref{hp:phi=}), (\ref{hp:phi-}), and (\ref{hp:phi3}).  Let us
	assume that $t_{0}>0$ satisfies (\ref{hp:t01}) and
	(\ref{hp:t02}).  Let $\beta(t)$, $\gamma(t)$, $\T$, $b(t)$,
	$c(t)$ be defined by (\ref{defn:beta-gamma}), (\ref{defn:T}),
	(\ref{defn:b}), and (\ref{defn:c}), respectively.
	
	Then there exists $u\in C^{2,1}(\T)$ satisfying equation
	\begin{equation}
		u_{t}(r,t)=-\varphi''(u_{r}(r,t))u_{rr}(r,t)-
		\frac{\varphi'(u_{r}(r,t))}{r}
		\label{eq:T-eq}
	\end{equation}
	for every $(r,t)\in\qtwo$, and the Neumann boundary conditions
	$$u_{r}(\beta(t_{0}-t),t)=u_{r}(\gamma(t_{0}-t),t)=1
		\quad\quad\quad
		 \forall t\in[0,t_{0}].$$
	
	Moreover we have that
	\begin{eqnarray*}
		 & u_{r}(r,t)>1
		\quad\quad\quad
		\forall (r,t)\in\intp(\T), &
		\label{th:T-sp}  \\
		 & u_{rr}(\beta(t_{0}-t),t)=b(t_{0}-t)
		\quad\quad\quad
		\forall t\in[0,t_{0}], &
		\label{th:T-urr-l}  \\
		 & u_{rr}(\gamma(t_{0}-t),t)=c(t_{0}-t)
		\quad\quad\quad
		\forall t\in[0,t_{0}]. &
		\label{th:T-urr-r}
	\end{eqnarray*}

\end{thm}

We didn't state explicitly what (\ref{g1:u}) and (\ref{g3:u}) become
in this new setting because, as we have seen in
Remark~\ref{rmk:heuristic}, they easily follow from the regularity of
the solution and the Neumann boundary conditions (\ref{g1:ur}) and
(\ref{g3:ur}).

In order to approach Theorem~\ref{thm:T}, for every $\ep\in(0,t_{0})$
we set (see Figure~\ref{fig:T})
\begin{equation}
	\Tep:=\left\{(r,t)\in\T:t\geq\ep\right\}.
	\label{defn:Tep}
\end{equation}

Now we approximate the degenerate parabolic problem in $\T$ with a
strictly parabolic problem in $\Tep$.  Note that in this case we 
prescribe an initial condition for $t=\ep$. This initial condition is 
compatible with the Neumann boundary conditions.

\begin{thm}\label{thm:Tep}
	Let $\varphi$, $t_{0}$, $\beta(t)$, $\gamma(t)$, $\T$,
	$b(t)$, $c(t)$ be as in Theorem~\ref{thm:T}. Let
	$\ep\in(0,t_{0})$, and let $\Tep$ be defined by (\ref{defn:Tep}).
	
	Then there exists a unique function $\uep\in C^{2,1}(\Tep)$
	satisfying equation (\ref{eq:T-eq}) for every $(r,t)\in\Tep$, the
	Neumann boundary conditions
	$$\uep_{r}(\beta(t_{0}-t),t)=\uep_{r}(\gamma(t_{0}-t),t)=
		1+\ep
		\quad\quad\quad
		 \forall t\in[\ep,t_{0}],$$
	and the initial condition
	$$\uep(r,\ep)=(1+\ep)r
		\quad\quad\quad
		 \forall r\in[\beta(t_{0}-\ep),\gamma(t_{0}-\ep)].$$

    Moreover $\uep$ satisfies the following estimates independent on
	$\ep$.
    \begin{enumerate}
		\item \emph{Maximum principle for space derivatives}.  We have that
		\begin{equation}
			1+\ep\leq\uep_{r}(r,t)\leq 3
			\quad\quad\quad
			\forall (r,t)\in\Tep.
			\label{est:Tep-ur}
		\end{equation}
		
		\item \emph{Uniform strict parabolicity in the interior}.  For
		every compact set $K\subseteq\Tep\cap\intp(\T)$ there exists a
		constant $M_{1}$ (depending on $K$ but not on $\ep$) such that
		\begin{equation}
			1<M_{1}\leq\uep_{r}(r,t)
			\hspace{3em}
			\forall(r,t)\in K.
			\label{est:Tep-unif-farab}
		\end{equation}
	
		As a consequence, there exists a constant $M_{2}>0$ such that
		\begin{equation}
			-\varphi''(\uep_{r}(r,t))\geq M_{2} \hspace{3em}
			\forall (r,t)\in K.
			\label{est:Tep-unif-parab}
		\end{equation}
	
		\item \emph{Estimate on second derivatives at the
		boundary}.  We have that
		\begin{equation}
	    	\left|\uep_{rr}(\beta(t_{0}-t),t)-b(t_{0}-t)\right|
			\leq \sqrt{\ep}
			\hspace{3em}
	        \forall t\in[\ep,t_{0}],
	    	\label{est:Tep-ueprr-l}
		\end{equation}
		\begin{equation}
	    	\left|\uep_{rr}(\gamma(t_{0}-t),t)-c(t_{0}-t)\right|
			\leq \sqrt{\ep}
			\hspace{3em}
	        \forall t\in[\ep,t_{0}].
	    	\label{est:Tep-ueprr-r}
		\end{equation}

		\item \emph{Global estimate on second derivatives}.  There
		exists a constant $M_{3}$ such that
		\begin{equation}
	    	\left|\uep_{rr}(r,t)-b(t_{0}-t)\right|
			\leq \sqrt{\ep}+M_{3}\left|r-\beta(t_{0}-t)\right|
			\hspace{3em}
	        \forall (r,t)\in\Tep,
	    	\label{est:Tep-ueprr-ln}
		\end{equation}
		\begin{equation}
	    	\left|\uep_{rr}(r,t)-c(t_{0}-t)\right|
			\leq \sqrt{\ep}+M_{3}\left|r-\gamma(t_{0}-t)\right|
			\hspace{3em}
	        \forall (r,t)\in\Tep.
	    	\label{est:Tep-ueprr-rn}
		\end{equation}
    
		As a consequence, there exist constants $M_{4}$ and $M_{5}$
		such that
		$$|\uep_{rr}(r,t)|\leq M_{4}
			\quad\quad\quad
			\forall (r,t)\in\Tep,$$
		$$|\uep_{t}(r,t)|\leq M_{5}
			\quad\quad\quad
			\forall (r,t)\in\Tep.$$
		
		\item \emph{Integral estimates}.  There exists a
		constant $M_{6}$ such that
		\begin{equation}
			-\int_{\Tep}
			\varphi''(\uep_{r}(r,t))
			[\uep_{rt}(r,t)]^{2}\,dr\,dt\leq M_{6}.
			\label{est:Tep-ueprt-int}
		\end{equation}
   
		In addition, for every closed rectangle
		$[r_{1},r_{2}]\times[t_{1},t_{2}]\subseteq\Tep\cap\intp(\T)$
		there exists a constant $M_{7}$ (depending on the rectangle,
		but not on $\ep$) such that
		$$\int_{r_{1}}^{r_{2}}[\uep_{rt}(r,t)]^{2}\,dr \leq M_{7}
			\hspace{3em} 
			\forall t\in[t_{1},t_{2}],$$
		$$\int_{r_{1}}^{r_{2}}[\uep_{rrr}(r,t)]^{2}\,dr \leq M_{7}
			\hspace{3em} 
			\forall t\in[t_{1},t_{2}],$$
		$$\int_{t_{1}}^{t_{2}}\int_{r_{1}}^{r_{2}}
	    	[\uep_{rrt}(r,t)]^{2}\,dr\,dt  \leq  M_{7}.$$
	\end{enumerate}
\end{thm}

\setcounter{equation}{0}
\section{Proof of estimates}\label{sec:proofs}

In the following two statements we collect the properties of $b(t)$
and $c(t)$ which are needed in the main proofs.  We only prove
Lemma~\ref{lemma:b} because the proof of Lemma~\ref{lemma:c} is quite
similar.

\begin{lemma}\label{lemma:b}
	Let $t_{0}>0$, and let $\beta(t)$ be defined by
	(\ref{defn:beta-gamma}). Let us assume that
	\begin{equation}
		\frac{1}{t_{0}}\geq 4
		\sqrt{\varphi'(1)\left|\varphi'''(1)\right|}.
		\label{hp:b}
	\end{equation}
	
	Then the function $b(t)$ introduced in (\ref{defn:b}) is well
	defined for every $t\in[0,t_{0}]$ and fulfils the following
	properties.
	\begin{enumerate}
		\renewcommand{\labelenumi}{(\arabic{enumi})}
		\item  We have that $b\in C^{0}([0,t_{0}])\cap
		C^{\infty}([0,t_{0}))$.
	
		\item  We have that
		\begin{equation}
			0<-b'(t)\leq 5t_{0}\varphi'(1)\beta'(t)
			\quad\quad
			\forall t\in[0,t_{0}).
			\label{est:b'}
		\end{equation}
	
		\item  The function $b(t)$ is decreasing and
		\begin{equation}
			0\leq b(t)\leq b(0)\leq 2\varphi'(1)t_{0}
			\quad\quad
			\forall t\in[0,t_{0}],
			\label{est:b-Q}
		\end{equation}
		\begin{equation}
			0\leq b(t_{0}-t)\leq 2\varphi'(1)\sqrt{t_{0}}\sqrt{t}
			\quad\quad
			\forall t\in[0,t_{0}].
			\label{est:b-T}
		\end{equation}
	\end{enumerate}
\end{lemma}

\prf
From (\ref{defn:beta-gamma}) we easily deduce that
\begin{equation}
	\beta(t)\geq 2,
	\hspace{3em}
	\beta'(t)\geq\frac{1}{2t_{0}},
	\hspace{3em}
	\beta''(t)=2t_{0}\left[\beta'(t)\right]^{3}
	\label{eq:prop-beta}
\end{equation}
for every $t\in[0,t_{0}]$. Moreover from assumption (\ref{hp:b}) we
have that
\begin{equation}
	\left[\beta'(t)\right]^{2}\geq 
	\frac{1}{4t_{0}^{2}}\geq 4\varphi'(1)|\varphi'''(1)|
	\quad\quad\quad
	\forall t\in[0,t_{0}),
	\label{delta}
\end{equation}

Let us consider the polynomial
$$p(x):=\varphi'''(1)x^{2}+\beta'(t)x-
\frac{\varphi'(1)}{\beta^{2}(t)},$$
where $t\in[0,t_{0}]$ is thought as a parameter.  Due to our
assumptions on $\varphi$ we have that $\varphi'''(1)\leq 0$.
Exploiting (\ref{delta}) one can therefore check that
\begin{equation}
	p\left(\frac{\varphi'(1)}{\beta^{2}(t)\beta'(t)}\right)\leq 0,
	\hspace{3em}
	p\left(\frac{\varphi'(1)}{\beta'(t)}\right)\geq 0.
	\label{est:p(x)}
\end{equation}

This implies that $p(x)$ has at least one real root.  Moreover $p(x)$
is either a polynomial of degree one, or a polynomial of degree two
representing a \emph{concave} function.  In both cases the unique or
the smallest root $b(t)$ is the one lying in the interval whose
endpoints appear in (\ref{est:p(x)}), hence
\begin{equation}
	\frac{\varphi'(1)}{\beta^{2}(t)\beta'(t)}\leq b(t)
	\leq\frac{\varphi'(1)}{\beta'(t)}
	\quad\quad
	\forall t\in[0,t_{0}).
	\label{th:est-b}
\end{equation}

The regularity of $b(t)$ in $[0,t_{0})$ follows from the regularity of
the coefficients of the equation.  The continuity up to $t=t_{0}$
follows from (\ref{th:est-b}) and the fact that $\beta'(t)\to +\infty$
as $t\to t_{0}^{-}$ (we remind that we set $b(t_{0})=0$). This proves 
statement (1).

Let us consider now the derivative $b'(t)$.  From the implicit
function theorem we have that
\begin{equation}
	-b'(t)=\frac{\beta''(t)b(t)+2\varphi'(1)\beta'(t)
	[\beta(t)]^{-3}}{2\varphi'''(1)b(t)+\beta'(t)}
	\quad\quad\quad
	\forall t\in[0,t_{0}).
	\label{b'}
\end{equation}

Let us estimate the numerator of (\ref{b'}).  Exploiting
(\ref{eq:prop-beta}) and the upper bound for $b(t)$ provided by
(\ref{th:est-b}), we have that
\begin{eqnarray*}
	0\ <\ \beta''(t)b(t)+2\frac{\varphi'(1)\beta'(t)}{\beta^{3}(t)}
	 & \leq & 2t_{0}[\beta'(t)]^{3}\frac{\varphi'(1)}{\beta'(t)}+
	 \frac{\varphi'(1)}{4}\beta'(t)
	  \\
	 & = & \varphi'(1)[\beta'(t)]^{2}\left(2t_{0}+
	 \frac{1}{4\beta'(t)}\right)
	  \\
	 & \leq & \frac{5t_{0}}{2}\varphi'(1)[\beta'(t)]^{2}.
\end{eqnarray*}

Let us estimate the denominator of (\ref{b'}).  Exploiting once again
the upper bound for $b(t)$ provided by (\ref{th:est-b}), and
inequality (\ref{delta}), we obtain that
$$2\varphi'''(1)b(t)+\beta'(t) \geq \beta'(t)-2
|\varphi'''(1)|\frac{\varphi'(1)}{\beta'(t)} = \beta'(t)\left(1-
\frac{2|\varphi'''(1)|\varphi'(1)}{[\beta'(t)]^{2}}\right)
\geq  \frac{1}{2}\beta'(t).$$

At this point (\ref{est:b'}) easily follows from the estimates on
the numerator and the denominator of (\ref{b'}). This completes the
proof of statement (2).

Let us prove now statement (3).  The function $b(t)$ is decreasing
because its derivative is negative.  It follows that $0= b(t_{0})\leq
b(t)\leq b(0)$ for every $t\in[0,t_{0}]$.  Finally, the estimates on
$b(0)$ and $b(t-t_{0})$ stated in (\ref{est:b-Q}) and (\ref{est:b-T})
follow from (\ref{eq:prop-beta}) and the upper bound in
(\ref{th:est-b}).  \qed

\begin{lemma}\label{lemma:c}
	Let $t_{0}$ be as in Lemma~\ref{lemma:b}, and let $\gamma(t)$ be
	defined by (\ref{defn:beta-gamma}).
	
	Then the function $c(t)$ introduced in (\ref{defn:c}) is well
	defined for every $t\in[0,t_{0}]$ and fulfils the following
	properties.
	\begin{enumerate}
		\renewcommand{\labelenumi}{(\arabic{enumi})}
		\item  We have that $c\in C^{0}([0,t_{0}])\cap
		C^{\infty}([0,t_{0}))$.
	
		\item  We have that (we remind that $\gamma'(t)$ is negative)
		$$0<c'(t)\leq -5t_{0}\varphi'(1)\gamma'(t)
		\quad\quad
		\forall t\in[0,t_{0}).$$
	
		\item  The function $c(t)$ is increasing and
		\begin{equation}
			-2\varphi'(1)t_{0}\leq c(0)\leq c(t)\leq 0 
			\quad\quad
			\forall t\in[0,t_{0}],
			\label{est:c-Q}
		\end{equation}
		\begin{equation}
			0\geq c(t_{0}-t)\geq -2\varphi'(1)\sqrt{t_{0}}\sqrt{t}
			\quad\quad
			\forall t\in[0,t_{0}].
			\label{est:c-T}
		\end{equation}
	\end{enumerate}
\end{lemma}

\begin{rmk}
	\begin{em}
		The second inequality in (\ref{hp:t01}) implies (\ref{hp:b}).
		We can therefore apply the conclusions of Lemma~\ref{lemma:b}
		and Lemma~\ref{lemma:c} in the proofs of our main results.
		Exploiting the first and second inequality in (\ref{hp:t01}),
		from statement~(3) in Lemma~\ref{lemma:b} and
		Lemma~\ref{lemma:c} we obtain the following stronger
		inequalities
		\begin{eqnarray}
			 & 0\leq b(t)\leq b(0)\leq 1
			 \quad\quad
			\forall t\in[0,t_{0}], & 
			\label{est:b-1}  \\
			 & 0\leq b(t_{0}-t)\leq \sqrt{t} 
			 \quad\quad
			\forall t\in[0,t_{0}],& 
			\label{est:b-sqrt}  \\
			 & -\sqrt{t}\leq c(t_{0}-t)\leq 0
			 \quad\quad
			\forall t\in[0,t_{0}]. & 
			\label{est:c-sqrt}
		\end{eqnarray}
	\end{em}
\end{rmk}

We finally state a classical comparison result for fully nonlinear
parabolic equations.  This is the key tool in our analysis.  We omit
the standard proof.

\begin{lemma}\label{lemma:super-sol}
	Let $\D\subseteq\re^{2}$ be a compact set, and let
	$\psi:\D\times\re^{3}\to\re$ be a continuous
	function. Let us assume that
	\begin{itemize}
		\item  $\psi$ is nondecreasing in the last variable
		(degenerate ellipticity), namely
		$$\psi(r,t,p,q,s_{1})\leq\psi(r,t,p,q,s_{2})
		\quad\quad
		\forall s_{1}\leq s_{2},\ 
		\forall(r,t,p,q)\in\D\times\re^{2},$$
	
		\item  $\psi$ is locally Lipschitz continuous in the third
		variable, namely for every $R\geq 0$ there exists a constant
		$L$ such that
		$$|\psi(r,t,p_{1},q,s)-\psi(r,t,p_{2},q,s)|\leq L|p_{1}-p_{2}|
		\quad\quad
		\forall(r,t,p_{1},p_{2},q,s)\in\D\times[0,R]^{4}.$$
	\end{itemize}
	
	Let $u$ and $v$ be two functions in $C^{0}(\D)\cap
	C^{2,1}(\intp(\D))$ such that
	$$u_{t}\leq\psi(r,t,u,u_{r},u_{rr})
	\hspace{3em}
	\forall(r,t)\in\intp(\D),$$
	$$v_{t}\geq\psi(r,t,v,v_{r},v_{rr})
	\hspace{3em}
	\forall(r,t)\in\intp(\D),$$
	$$u(r,t)\leq v(r,t)
	\hspace{3em}
	\forall(r,t)\in\partial_{P}(\D).$$
        
    Then $u(r,t)\leq v(r,t)$ for every $(r,t)\in\D$.
	\qed
\end{lemma}

\begin{rmk}\label{rmk:super-sol}
	\begin{em}
		If we know a priori that there exist constants $a$ and $b$
		such that $a\leq u(r,t)\leq b$ for every $(r,t)\in\D$, then we
		can weaken the assumptions of Lemma~\ref{lemma:super-sol} by
		asking that $v$ fulfils its differential inequality only for
		those $(r,t)\in\intp(\D)$ such that $a\leq v(r,t)\leq b$.  
 		We can obtain a similar statement by swapping the role
 		of $u$ and $v$.
	\end{em}
\end{rmk}

\subsection{Proof of Theorem~\ref{thm:Qep}}

Let us briefly sketch the outline of the proof.  First of all we show
that the initial boundary value problem has a unique solution $\uep\in
C^{2,1}(\qone)\cap C^{\infty}(\intp(\qone))$.  Then for the sake of
simplicity we set $v=:\uep_{r}$, and $w:=\uep_{rr}$.  These functions
belong to $C^{0}(\qone)\cap C^{\infty}(\intp(\qone))$.  It is easy
to see that $v$ is the solution of the following equation
\begin{equation}
	v_{t}=\varphi''(v)v_{rr}+\varphi'''(v)v_{r}^{2}+
	\frac{\varphi''(v)}{r}v_{r}-\frac{\varphi'(v)}{r^{2}}
	\quad\quad\quad
	\forall (r,t)\in\intp(\qone),
	\label{eq:v}
\end{equation}
with Dirichlet boundary conditions
\begin{equation}
	v(1,t)=0
	\quad\quad\quad
	\forall t\in[0,t_{0}],
	\label{eq:v-bcl}
\end{equation}
\begin{equation}
	v(\beta(t),t)=1-\ep
	\quad\quad\quad
	\forall t\in[0,t_{0}],
	\label{eq:v-bcr}
\end{equation}
and initial datum 
\begin{equation}
	v(r,0)=(1-\ep)u_{0r}(r)
	\quad\quad\quad
	\forall r\in[1,2].
	\label{eq:v-cauchy}
\end{equation}

In the same way $w$ is a solution in $\intp(\qone)$ of equation
\begin{eqnarray}
	w_{t} & = &
	\varphi''(v)w_{rr}+3\varphi'''(v)w_{r}w+\varphi^{IV}(v)w^{3}
	\nonumber  \\
	 & & + \frac{\varphi'''(v)}{r}w^{2}+\frac{\varphi''(v)}{r}w_{r}
	 -2\frac{\varphi''(v)}{r^{2}}w+2\frac{\varphi'(v)}{r^{3}},
	\label{eq:w}
\end{eqnarray}
where the terms in $v$ are to be interpreted as coefficients depending
on $r$ and $t$. Moreover $w$ satisfies the Dirichlet boundary conditions
\begin{equation}
	w(1,t)=\uep_{rr}(1,t)
	\quad\quad
	\forall t\in[0,t_{0}],
	\label{eq:w:bc1}
\end{equation}
\begin{equation}
	w(\beta(t),t)=\uep_{rr}(\beta(t),t)
	\quad\quad
	\forall t\in[0,t_{0}],
	\label{eq:w:bc2}
\end{equation}
and the initial condition
\begin{equation}
	w(r,0)=(1-\ep)u_{0rr}(r)
	\quad\quad
	\forall r\in[1,2].
	\label{eq:w-cauchy}
\end{equation}

The initial data (\ref{eq:v-cauchy}) and (\ref{eq:w-cauchy}) can
be easily estimated using (\ref{hp:u0-comp}) through (\ref{hp:u0rrr}).
Indeed from Taylor's expansion we have that
$$u_{0r}(r)=u_{0r}(2)+u_{0rr}(2)(r-2)+\frac{1}{2}u_{0rrr}(\xi)
(r-2)^{2}$$
for a suitable $\xi\in(1,2)$, hence
\begin{equation}
	1+b(0)(r-2)-5(r-2)^{2}\leq u_{0r}(r)\leq
	1+b(0)(r-2)+5(r-2)^{2} 
	\quad\quad 
	\forall r\in[1,2].
	\label{est:u0r}
\end{equation}

Analogously we have that
$u_{0rr}(r)=u_{0rr}(2)+u_{0rrr}(\xi)(r-2)$,
hence
\begin{equation}
	b(0)-10(2-r)\leq u_{0rr}(r)\leq
	b(0)+10(2-r) 
	\quad\quad\quad 
	\forall r\in[1,2].
	\label{est:u0rr}
\end{equation}

Then we show that $v$ and $w$ satisfy five sets of inequalities in
$\qone$. 
\begin{itemize}
	\item The first pair of inequalities is
	\begin{equation}
		0\leq v(r,t)\leq 1-\ep,
		\label{est:ur-mp}
	\end{equation}
	which is equivalent to (\ref{est:Q-ur}).

	\item The second inequality is
	\begin{equation}
		v(r,t)\leq 1-\eta\left((r-3)^{2}+\frac{t}{t_{0}}-1\right)
		\label{est:Qep-up}
	\end{equation}
	for a suitable constants $\eta>0$.  The term after $\eta$ is
	positive in $\qone\setminus\Gamma_{1}$.  This implies
	(\ref{est:Q-unif-farab}), hence also (\ref{est:Q-unif-parab}).

	\item The third pair of inequalities is
	\begin{equation}
		0\leq v(r,t)\leq \frac{20}{(1-800\kf t)^{1/2}}
		\left(1-e^{1-r}\right),
		\label{est:ueprr-l}
	\end{equation}
	where $\gamma_{2}$ is the constant defined in (\ref{defn:gamma2}).
	The left-hand side and the right-hand side coincide for $r=1$,
	hence $v_{r}(1,t)$ is bounded by their derivative computed in
	$r=1$.  Thus we obtain that $0\leq v_{r}(1,t)\leq 20(1-800\kf
	t)^{-1/2}$.  The upper bound is less than 100 
	due to the second inequality in (\ref{hp:t01}).  This proves
	(\ref{est:Q-ueprr-l}).

	\item The fourth pair of inequalities, and probably the most
	delicate one, is
	\begin{equation}
		v(r,t)\geq
		1-\ep-(b(t)+\ep)(\beta(t)-r)-\gamma_{0}(\beta(t)-r)^{2},
		\label{est:ueprr-rl}
	\end{equation}
	\begin{equation}
		v(r,t)\leq
		1-\ep-(b(t)-\ep)(\beta(t)-r)+\gamma_{0}(\beta(t)-r)^{2},
		\label{est:ueprr-ru}
	\end{equation}
	where $\gamma_{0}$ is the constant defined in (\ref{defn:gamma0}).
	In this case $v(r,t)$ is bounded from below and from above by two
	functions which coincide for $r=\beta(t)$.  It follows that
	$v_{r}(\beta(t),t)$ is bounded by the space derivatives of these
	two functions computed in $r=\beta(t)$.  This implies
	(\ref{est:Q-ueprr}). 

	\item Thanks to (\ref{est:Q-ueprr-l}) and (\ref{est:Q-ueprr}) we
	have an estimate on the values of $w$ at the boundary.  This is
	the starting point to prove the fifth set of inequalities
	\begin{equation}
		b(t)-\ep-\gamma_{1}(\beta(t)-r)\leq w(r,t)\leq
		b(t)+\ep+\gamma_{1}(\beta(t)-r)
		\label{est:ueprr}
	\end{equation}
	where $\gamma_{1}$ is the constant defined in (\ref{defn:gamma1}).
	This implies (\ref{est:Q-ueprr-n}), hence also
	(\ref{est:ueprr-glob}) and (\ref{est:uept-glob}).
\end{itemize}

All these estimates are proved using subsolutions and supersolutions.
Finally the proof of (\ref{est:ueprt-int}) through
(\ref{est:ueprrt-int}) relies on usual energy estimates.

We are now ready to proceed with the details.

\paragraph{Existence and maximum principle for space derivatives}

Let us take a function $\varphi_{\ep}\in C^{\infty}(\re)$ which
coincides with $\varphi$ in the interval $[0,1-\ep]$, and such that
$\varphi_{\ep}''(\sigma)\geq\nu_{\ep}>0$ for every $\sigma\in\re$.
Equation (\ref{eq:PM-rad}), with $\varphi_{\ep}$ instead of $\varphi$,
is strictly parabolic.  The initial boundary value problem in $\qone$
can be reduced to the fixed domain $[0,1]\times[0,t_{0}]$ by the
variable change $(s,t)\to(1+(\beta(t)-1)s,t)$, which of course doesn't
change the strict parabolicity.  By well knows results (see for
example \cite{lady}) the problem admits a unique solution $\uep\in
C^{2,1}(\qone)\cap C^{\infty}(\intp(\qone))$.  By the way, the
estimates we are going to prove could be used to give a self contained
proof of the existence result for $\ep$ fixed, showing in particular
that the singularity of $\qone$ in $(3,t_{0})$ doesn't affect the
existence or the regularity of the solution.

If we show that this solution satisfies (\ref{est:Q-ur}), then this
same solution satisfies equation (\ref{eq:PM-rad}) with the original
$\varphi$.  Our assumptions on $\varphi$ imply that $\varphi'(0)=0$
and $\varphi'(1-\ep)>0$.  It follows that $z(r,t):=0$ is a subsolution
of (\ref{eq:v}) through (\ref{eq:v-cauchy}), while $z(r,t):=1-\ep$ is
a supersolution of the same problem.  The usual comparison principle
implies (\ref{est:ur-mp}).

\paragraph{Uniform strict parabolicity in the interior}

Let us choose $\eta>0$ small enough so that
\begin{equation}
	\eta\leq \frac{1}{8},
	\hspace{2em}
	\eta\leq \inf_{r\in[1,2)}\frac{1-u_{0r}(r)}{(2-r)(4-r)},
	\hspace{3em}
	\eta\leq\frac{1}{9}\left(\frac{1}{t_{0}}+20\gamma_{2}\right)^{-1}
	\cdot\varphi'\left(\frac{1}{2}\right).
	\label{defn:eta}
\end{equation}

Note that the infimum is positive because our assumptions on $u_{0}$
imply that the fraction is positive for every $r\in[1,2)$ and tends to
$b(0)/2>0$ as $r\to 2^{-}$.

Let $z(r,t)$ denote the right-hand side of (\ref{est:Qep-up}).  We
claim that, for this choice of $\eta$, $z$ is a supersolution
of (\ref{eq:v}) through (\ref{eq:v-cauchy}), which implies
(\ref{est:Qep-up}).

\subparagraph{\emph{\textmd{Boundary $r=1$}}}

Since $|r-3|\leq 2$, $t\leq t_{0}$, and $\eta\leq 1/8$, we have that
\begin{equation}
	z(r,t)=1-\eta\left((r-3)^{2}+\frac{t}{t_{0}}-1\right)\geq
	1-4\eta\geq \frac{1}{2}
	\quad\quad
	\forall (r,t)\in\qone.
	\label{est:up-z-glob}
\end{equation}

In particular we have that $z(1,t)\geq 0=v(1,t)$.

\subparagraph{\emph{\textmd{Boundary $r=\beta(t)$}}}

In this case we have that $z(\beta(t),t)=1\geq 1-\ep=v(\beta(t),t)$.

\subparagraph{\emph{\textmd{Boundary $t=0$}}}

Thanks to the second inequality in (\ref{defn:eta}), and the fact that
$u_{0r}(r)\geq 0$ for every $r\in[1,2]$, we have that 
$$z(r,0)=1-\eta\left((r-3)^{2}-1\right)= 1-\eta(2-r)(4-r)\geq$$
$$\geq u_{0r}(r)\geq
(1-\ep)u_{0r}(r)=v(r,0).$$

\subparagraph{\emph{\textmd{Differential inequality}}}

We have that
$$z_{t}=-\frac{\eta}{t_{0}},
\hspace{2em}
z_{r}(r,t)=-2\eta(r-3),
\hspace{2em}
z_{rr}(r,t)=-2\eta.$$

Therefore we have to prove that
$$-\frac{\eta}{t_{0}}\geq -2\eta\varphi''(z)+
\varphi'''(z)4\eta^{2}(r-3)^{2}- \frac{\varphi''(z)}{r}2\eta(r-3)
-\frac{\varphi'(z)}{r^{2}}.$$

Due to (\ref{est:up-z-glob}) and the properties of $\varphi'$, we have
that $\varphi'(z)/r^{2}\geq (1/9)\varphi'(1/2)$.  All the other terms
are uniformly small when $\eta$ is small.  This shows that the
differential inequality is satisfied when $\eta$ is small enough, for
example as soon as $\eta$ fulfils the last inequality in
(\ref{defn:eta}).

\paragraph{Estimate on second derivatives at the fixed boundary} 

Let us prove (\ref{est:ueprr-l}).  We already know from
(\ref{est:ur-mp}) that $v(r,t)\geq 0$ in $\qone$.  In order to prove
the other inequality, we set for simplicity $k(t):=20(1-800\kf t
)^{-1/2}$.  It is easy to check that $k(t)$ is the solution of the
Cauchy problem 
$$k'(t)=\kf k^{3}(t), \quad\quad\quad
k(0)=20.$$

We claim that $z(r,t):=k(t)\left(1-e^{1-r}\right)$ is a supersolution
of (\ref{eq:v}) through (\ref{eq:v-cauchy}), which implies
(\ref{est:ueprr-l}).

\subparagraph{\emph{\textmd{Boundary $r=1$}}}

In this case we have that $z(1,t)=0=v(1,t)$.

\subparagraph{\emph{\textmd{Boundary $r=\beta(t)$}}}

Since $\beta(t)\geq 2$, in this case we have that 
$$z(\beta(t),t)=k(t)\left(1-e^{1-\beta(t)}\right)
\geq 20(1-e^{-1})\geq 1\geq v(\beta(t),t).$$   

\subparagraph{\emph{\textmd{Boundary $t=0$}}}

From the assumptions on $u_{0}$ we have that
\begin{equation}
	v(r,0)=(1-\ep)u_{0r}(r)\leq u_{0r}(r)=
	u_{0r}(1)+(r-1)u_{0rr}(\xi)\leq 10(r-1).
	\label{est:ur-vr0}
\end{equation}

On the other hand, the function $1-e^{1-r}$ is concave, hence
\begin{equation}
	z(r,0)=20\left(1-e^{1-r}\right)\geq
	20(1-e^{-1})(r-1)\geq 10(r-1).
	\label{est:ur-zr0}
\end{equation}

From (\ref{est:ur-vr0}) and (\ref{est:ur-zr0}) it follows that
$z(r,0)\geq v(r,0)$ for every $r\in[1,2]$.

\subparagraph{\emph{\textmd{Differential inequality}}}

We have that
$$z_{t}(r,t)=k'(t)\left(1-e^{1-r}\right),
\hspace{2em}
z_{r}(r,t)=k(t)e^{1-r},
\hspace{2em}
z_{rr}(r,t)=-k(t)e^{1-r}.$$

We have therefore to prove that
\begin{equation}
	k'(t)\left(1-e^{1-r}\right)\geq \varphi'''(z)k^{2}(t)e^{2-2r}
	-k(t)\left(1-\frac{1}{r}\right)\varphi''(z)e^{1-r}
	-\frac{\varphi'(z)}{r^{2}}.
	\label{eq:ur-di}
\end{equation}

Thanks to Remark~\ref{rmk:super-sol}, this inequality has to be
satisfied only for those values
$(z,t)\in\qone$ for which $0\leq z(r,t)\leq 1$. In this case we know
that the last two terms are negative. Since $\varphi'''(0)=0$ we have also
that
$$\left|\varphi'''(z)\right|=\left|\varphi'''(z)-\varphi'''(0)\right|
=\left|\varphi^{IV}(\xi)\right|z\leq
\kf k(t)\left(1-e^{1-r}\right).$$

It follows that
$$\mbox{right-hand side of (\ref{eq:ur-di})}\leq
\left|\varphi'''(z)\right|
k^{2}(t)\leq \gamma_{2}k^{3}(t)\left(1-e^{1-r}\right)=
k'(t)\left(1-e^{1-r}\right),$$
which completes the proof of the differential inequality.

\paragraph{Estimate on second derivatives at the moving boundary --
Subsolution} 

Let us prove (\ref{est:ueprr-rl}).  To this end, we denote its
right-hand side by $z(r,t)$, and we show that $z$ is a subsolution of
(\ref{eq:v}) through (\ref{eq:v-cauchy}).

\subparagraph{\emph{\textmd{Boundary $r=1$}}}

Since $\beta(t)\geq 2$ and $\gamma_{0}\geq 1$, in this case we have
that
$$z(1,t)=1-\ep-(b(t)+\ep)(\beta(t)-1)-
\gamma_{0}(\beta(t)-1)^{2}\leq 1-\gamma_{0}\leq 0=v(1,t).$$

\subparagraph{\emph{\textmd{Boundary $r=\beta(t)$}}}

In this case we have that $z(\beta(t),t)=1-\ep=v(\beta(t),t)$.

\subparagraph{\emph{\textmd{Boundary $t=0$}}}

Since $b(0)\geq 0$ and $\gamma_{0}\geq 5$, from (\ref{est:u0r}) have
that 
$$v(r,0) = (1-\ep)u_{0r}(r) \geq
(1-\ep)\left(1-b(0)(2-r)-5(2-r)^{2}\right) \geq$$
$$\geq 1-\ep-(b(0)+\ep)(2-r)-\gamma_{0}(2-r)^{2} = z(r,0).$$

\subparagraph{\emph{\textmd{Differential inequality}}}

Let us set for simplicity $x:=\beta(t)-r$. Then 
we have that
$$z_{r}(r,t)=b(t)+\ep+2\gamma_{0}x,
\hspace{3em}
z_{rr}(r,t)=-2\gamma_{0},$$
$$z_{t}(r,t)=-b(t)\beta'(t)-\beta'(t)\ep-
(b'(t)+2\gamma_{0}\beta'(t))x.$$

The differential inequality is therefore the following
\begin{eqnarray*}
	-b(t)\beta'(t)-\beta'(t)\ep-(b'(t)+2\gamma_{0}\beta'(t))x & \leq
	&
	-2\gamma_{0}\varphi''(z)+\varphi'''(z)(b(t)+\ep+2\gamma_{0}x)^{2}
	\\
	 &  & +
	\frac{\varphi''(z)}{r}(b(t)+\ep+2\gamma_{0}x)-\frac{\varphi'(z)}{r^{2}}.
\end{eqnarray*}

Thanks to Remark~\ref{rmk:super-sol} we can limit ourselves to prove
it for all values $(z,t)\in\qone$ such that $0\leq z(r,t)\leq 1$.
From the definition of $b(t)$ we have that
\begin{equation}
	-b(t)\beta'(t)=\varphi'''(1)b^{2}(t)-\frac{\varphi'(1)}{\beta^{2}(t)}.
	\label{eq-beta-appl}
\end{equation}

After changing the signs we have therefore to prove that
\begin{eqnarray}
	\beta'(t)\ep+(b'(t)+2\gamma_{0}\beta'(t))x & \geq & 
	2\gamma_{0}\varphi''(z)
	\nonumber  \\
	 &  & -\left[\varphi'''(z)(b(t)+\ep +2\gamma_{0}x)^{2}-
	 \varphi'''(1)b^{2}(t)\right]  
	\nonumber  \\
	 &  & -\frac{\varphi''(z)}{r}(b(t)+\ep +2\gamma_{0}x)  
	\nonumber  \\
	 &  & +\left[\frac{\varphi'(z)}{r^{2}}
	    -\frac{\varphi'(1)}{\beta^{2}(t)}\right]
	\nonumber  \\
	 & =: & I_{1}+I_{2}+I_{3}+I_{4}.
	\label{defn:urr-i}
\end{eqnarray}

Let us estimate the left-hand side.  Exploiting (\ref{est:b'}), and
the fact that $t_{0}\leq 1$ and $2\gamma_{0}\geq 2+5\varphi'(1)$, we
have that
\begin{eqnarray}
	\beta'(t)\ep+(b'(t)+2\gamma_{0}\beta'(t))x & \geq &
	\beta'(t)\ep+\left(-5\varphi'(1)t_{0}+2\gamma_{0}\right)\beta'(t)x
	\nonumber  \\
	 & \geq & \beta'(t)\ep+	2\beta'(t)x
	\nonumber  \\
	 & \geq & \frac{\ep}{2t_{0}}+\frac{x}{t_{0}}.
	\label{est:urr-lhs}
\end{eqnarray}

Let us estimate the four terms in the right-hand side of
(\ref{defn:urr-i}). From now on we exploit several times 
the following inequalities
$$0<\ep< 1,
	\hspace{2em}
	0\leq x\leq 2,
	\hspace{2em}
	0\leq z\leq 1,
	\hspace{2em}
	0\leq b(t)\leq 1,
	\hspace{2em}
	1\leq r\leq\beta(t).$$

First of all we remark that
$$|z-1|\leq\ep+(b(t)+\ep+\gamma_{0}x)x\leq\ep+(2\gamma_{0}+2)x.$$

As a consequence we have that
\begin{equation}
	|\varphi''(z)|=|\varphi''(z)-\varphi''(1)|=
	|z-1|\cdot|\varphi'''(\xi)|\leq\kf\ep+\kf(2\gamma_{0}+2)x,
	\label{est:varphi''}
\end{equation}
hence
\begin{equation}
	I_{1}\leq 2\gamma_{0}|\varphi''(z)|\leq 2\kf\gamma_{0}\ep+
	\kf(4\gamma_{0}^{2}+4\gamma_{0})x.
	\label{est:urr-i1}
\end{equation}

Analogously we have that
$$\left|\varphi'''(z)-\varphi'''(1)\right|=|z-1|\cdot
|\varphi^{IV}(\xi)|\leq
\kf\ep+\kf(2\gamma_{0}+2)x,$$
and
\begin{eqnarray*}
	\left|(b(t)+\ep+2\gamma_{0}x)^{2}-b^{2}(t)\right| & = &
	(2b(t)+\ep)\ep+ (4\gamma_{0}^{2}x+4b(t)\gamma_{0}+4\ep\gamma_{0})x \\
	 & \leq & 3\ep+(8\gamma_{0}^{2}+8\gamma_{0})x,
\end{eqnarray*}
hence
\begin{eqnarray}
	I_{2} & \leq & \left|\varphi'''(z)\right|\cdot
	\left|(b(t)+\ep +2\gamma_{0}x)^{2}-b^{2}(t)\right|+
	b^{2}(t)\left|\varphi'''(z)-\varphi'''(1)\right|
	\nonumber  \\
	 & \leq & 4\kf\ep+\kf(8\gamma_{0}^{2}+10\gamma_{0}+2)x.
	\label{est:urr-i2}
\end{eqnarray}

Exploiting (\ref{est:varphi''}) once again, we obtain that
\begin{eqnarray}
	I_{3} & \leq & (b(t)+\ep+2\gamma_{0}x)\frac{|\varphi''(z)|}{r}
	\nonumber  \\
	 & \leq & (4\gamma_{0}+2)|\varphi''(z)|
	\nonumber  \\
	 & \leq & \kf(4\gamma_{0}+2)\ep+\kf(8\gamma_{0}^{2}+
	 12\gamma_{0}+4)x.
	\label{est:urr-i3}
\end{eqnarray}

Finally we have that
$$|\varphi'(z)-\varphi'(1)|=
|z-1|\cdot|\varphi''(\xi)|\leq\kf\ep+\kf(2\gamma_{0}+2)x,$$
and
$$\left|\frac{1}{r^{2}}-\frac{1}{\beta^{2}(t)}\right|=
\frac{|\beta(t)+r|\cdot|\beta(t)-r|}{r^{2}\beta^{2}(t)}\leq
\frac{|\beta(t)+r|}{\beta^{2}(t)}x\leq x,$$
hence
\begin{eqnarray}
	I_{4} & \leq & |\varphi'(z)|\cdot\left|\frac{1}{r^{2}}-
	\frac{1}{\beta^{2}(t)}\right|+\frac{1}{\beta^{2}(t)}
	\left|\varphi'(z)-\varphi'(1)\right|
	\nonumber  \\
	 & \leq & \kf\ep+\kf(2\gamma_{0}+3)x.
	\label{est:urr-i4}  
\end{eqnarray}

Summing (\ref{est:urr-i1}) through (\ref{est:urr-i4}), and exploiting
(\ref{hp:t02}), we obtain that
\begin{equation}
	I_{1}+I_{2}+I_{3}+I_{4}\leq \kf(6\gamma_{0}+7)\ep+
	\kf(20\gamma_{0}^{2}+28\gamma_{0}+9)x\leq
	\frac{\ep}{2t_{0}}+\frac{x}{t_{0}}.
	\label{est:urr-rhs}
\end{equation}

This estimate and (\ref{est:urr-lhs}) imply the differential
inequality.

\paragraph{Estimate on second derivatives at the moving boundary --
Supersolution} 

Let us prove (\ref{est:ueprr-ru}).  To this end we show that its
right-hand side, which we denote by $z(r,t)$, is a supersolution of
(\ref{eq:v}) through (\ref{eq:v-cauchy}).

\subparagraph{\emph{\textmd{Boundary $r=1$}}}

Since $\gamma_{0}\geq 1\geq b(t)$ and $\beta(t)\geq 2$, in
this case we have that 
$$z(1,t)=1-\ep+(\beta(t)-1)
\left[\gamma_{0}(\beta(t)-1)-b(t)+\ep\right]\geq 1-\ep=v(1,t).$$

\subparagraph{\emph{\textmd{Boundary $r=\beta(t)$}}}

As in the case of the subsolution we have that
$$z(\beta(t),t)=1-\ep=v(\beta(t),t).$$

\subparagraph{\emph{\textmd{Boundary $t=0$}}}

Since $\gamma_{0}\geq 5$ and $b(0)\leq 1$, from (\ref{est:u0r}) we
have that
$$v(r,0)=(1-\ep)u_{0r}(r)\leq
(1-\ep)\left[1+b(0)(r-2)+5(r-2)^{2}\right]\leq$$
$$\leq 1-\ep+(b(0)-\ep)(r-2)+\gamma_{0}(r-2)^{2} =z(r,0).$$

\subparagraph{\emph{\textmd{Differential inequality}}}

We limit ourselves to sketch the argument, which is analogous to the
case of the subsolution.  After computing the derivatives and using
equation (\ref{eq-beta-appl}), we reduce ourselves to prove that (as
before $x:=\beta(t)-r$)
\begin{eqnarray*}
	\beta'(t)\ep+(2\gamma_{0}\beta'(t)-b'(t))x & \geq & 
	2\gamma_{0}\varphi''(z)
	\nonumber  \\
	 &  & +\left[\varphi'''(z)(b(t)-\ep-2\gamma_{0}x)^{2}-
	 \varphi'''(1)b^{2}(t)  \right]
	\nonumber  \\
	 &  & +\frac{\varphi''(z)}{r}(b(t)-\ep-2\gamma_{0}x)  
	\nonumber  \\
	 &  & +\left[\frac{\varphi'(1)}{\beta^{2}(t)}-
	 \frac{\varphi'(z)}{r^{2}}  \right]
	\nonumber  \\
	 & =: & I_{1}+I_{2}+I_{3}+I_{4}.
	\label{defn:urr+-i}
\end{eqnarray*}

Let us estimate the left-hand side.  Since $b'(t)\leq 0$ and
$\gamma_{0}\geq 1$, we have that
$$\beta'(t)\ep+(2\gamma_{0}\beta'(t)-b'(t))x\geq
\beta'(t)\ep+2\beta'(t)x\geq \frac{\ep}{2t_{0}}+\frac{x}{t_{0}}.$$

The four terms in the right-hand side can be estimated as in the case
of the subsolution.  In this way we obtain (\ref{est:urr-rhs}) once
again.  We conclude by using the smallness of $t_{0}$ exactly as in
the case of the subsolution.

\paragraph{Global estimate on second derivatives -- Supersolution} 

Let us prove the upper bound in (\ref{est:ueprr}). To this end we
denote the right-hand side by $z(r,t)$, and we show that it is
a supersolution of (\ref{eq:w}) through (\ref{eq:w-cauchy}).

\subparagraph{\emph{\textmd{Boundary $r=1$}}}

Since $b(t)\geq 0$, $\beta(t)\geq 2$, and $\gamma_{1}\geq 100$, from
(\ref{est:Q-ueprr-l}) we have that
$$z(1,t)=b(t)+\ep+\gamma_{1}(\beta(t)-1)\geq \gamma_{1}\geq
100\geq\uep_{rr}(1,t)=w(1,t).$$

\subparagraph{\emph{\textmd{Boundary $r=\beta(t)$}}}

From (\ref{est:Q-ueprr}) we have that
$$z(\beta(t),t)=b(t)+ \ep\geq \uep_{rr}(\beta(t),t)=
w(\beta(t),t).$$

\subparagraph{\emph{\textmd{Boundary $t=0$}}}

Since $\gamma_{1}\geq 10$, from (\ref{est:u0rr}) we have that
$$z(r,0)\geq b(0)+\gamma_{1}(2-r)\geq b(0)+10(2-r)\geq
(1-\ep)\left[b(0)+10(2-r)\right]\geq$$
$$\geq(1-\ep)u_{0rr}(r)=w(r,0).$$

\subparagraph{\emph{\textmd{Differential inequality}}}

We have that 
$$z_{r}(r,t)=-\gamma_{1},
\hspace{2em} 
z_{rr}(r,t)=0, 
\hspace{2em}
z_{t}(r,t)=b'(t)+\gamma_{1}\beta'(t).$$

The inequality to be satisfied in $\qone$ is therefore the following
$$b'(t)+\gamma_{1}\beta'(t)\geq -3\gamma_{1}\varphi'''(v)z+
\varphi^{IV}(v)z^{3}+\frac{\varphi'''(v)}{r}z^{2}-
\gamma_{1}\frac{\varphi''(v)}{r}
-2\frac{\varphi''(v)}{r^{2}}z+2\frac{\varphi'(v)}{r^{3}},$$
where of course the terms in $v$ are thought as coefficients.

Let us examine the left-hand side.  Using (\ref{est:b'}) and the fact
that $t_{0}\leq 1$ and $\gamma_{1}\geq 5\varphi'(1)+1$, we have that
\begin{equation}
	b'(t)+\gamma_{1}\beta'(t)\geq
	\left(-5\varphi'(1)t_{0}+\gamma_{1}\right)
	\beta'(t)\geq\beta'(t)\geq\frac{1}{2t_{0}}.
	\label{est:w-lhs}
\end{equation}

Since $|z(r,t)|\leq 2+2\gamma_{1}$, a rough estimate of the right-hand side 
of the differential inequality gives that
$$\mbox{right-hand side}\leq \kf(4\gamma_{1}+6)(2+2\gamma_{1})^{3}
\leq 48\kf(\gamma_{1}+1)^{4}.$$

Due to the last condition in (\ref{hp:t01}), this estimate and
(\ref{est:w-lhs}) imply the differential inequality.

\paragraph{Global estimate on second derivatives -- Subsolution} 

Let us prove the lower bound in (\ref{est:ueprr}). To this end we
denote the left-hand side by $z(r,t)$, and we show that it is
a subsolution of (\ref{eq:w}) through (\ref{eq:w-cauchy}).

\subparagraph{\emph{\textmd{Boundary $r=1$}}}

Since $\beta(t)-1\geq 1$ and $b(t)\leq 1\leq\gamma_{1}$, from
(\ref{est:Q-ueprr-l}) we have that
$$z(1,t)=b(t)-\ep-\gamma_{1}(\beta(t)-1)\leq 1-\gamma_{1}\leq 0\leq
\uep_{rr}(1,t)=w(1,t).$$

\subparagraph{\emph{\textmd{Boundary $r=\beta(t)$}}}

From (\ref{est:Q-ueprr}) we have that
$$z(\beta(t),t)=b(t)-\ep \leq \uep_{rr}(\beta(t),t)=
w(\beta(t),t).$$

\subparagraph{\emph{\textmd{Boundary $t=0$}}}

Since $b(0)\leq 1$ and $\gamma_{1}\geq 10$, from (\ref{est:u0rr}) we
have that
$$z(r,0)= b(0)-\ep-\gamma_{1}(2-r)\leq 
(1-\ep)(b(0)-10(2-r))\leq
(1-\ep)u_{0rr}(r)=w(r,0).$$

\subparagraph{\emph{\textmd{Differential inequality}}}

We have to prove that 
$$b'(t)-\gamma_{1}\beta'(t)\leq
3\gamma_{1}\varphi'''(v)z+
\varphi^{IV}(v)z^{3}+\frac{\varphi'''(v)}{r}z^{2}+
\gamma_{1}\frac{\varphi''(v)}{r}
-2\frac{\varphi''(v)}{r^{2}}z+2\frac{\varphi'(v)}{r^{3}}.$$

Since $\gamma_{1}\geq 1$ and $b'(t)\leq 0$, in the left-hand side we
have that
$$b'(t)-\gamma_{1}\beta'(t)\leq -\beta'(t)\leq-\frac{1}{2t_{0}},$$
while the usual rough estimate on the right-hand side gives that
$$\mbox{right-hand side} \geq -48\kf(\gamma_{1}+1)^{4}.$$

The conclusion follows as in the preceding case.

\paragraph{Integral estimates}

In the following estimates we introduce constants $c_{1}$, $c_{2}$, \ldots,
all independent on $\ep$.

Computing the time derivative and integrating by parts, we have that
\begin{eqnarray}
    \frac{\mbox{d}}{\mbox{d}t}\left(
    \int_{1}^{\beta(t)}\left[\uep_{t}\right]^{2}
    dr\right) & = & 2\int_{1}^{\beta(t)}\uep_{t}\cdot
    \frac{1}{r}\left(r\varphi'(\uep_{r})\right)_{rt}dr+
    \beta'(t)\left[\uep_{t}(\beta(t),t)\right]^{2}
    \nonumber  \\
     & = & -2\int_{1}^{\beta(t)}\varphi''(\uep_{r})
     \left[\uep_{rt}\right]^{2}\,dr+
     2\int_{1}^{\beta(t)}\frac{\varphi''(\uep_{r})}{r}
     \uep_{t}\uep_{rt}\,dr
    \nonumber  \\
     &  & +2\varphi''(1-\ep)\uep_{t}(\beta(t),t)
     \uep_{rt}(\beta(t),t)+
    \beta'(t)\left[\uep_{t}(\beta(t),t)\right]^{2}
    \nonumber  \\
     & =: & I_{1}+I_{2}+I_{3}+I_{4}.
    \label{est:int-i}
\end{eqnarray}

Note that in the integration by parts we neglected the boundary term
in $r=1$ because $\uep_{rt}(1,t)=0$ due to the Neumann boundary
condition (\ref{eq:nbc-l}).

Let us estimate some of the terms in (\ref{est:int-i}). From
(\ref{est:uept-glob}) we have that
$$2\frac{\varphi''(\uep_{r})}{r}\uep_{t}\uep_{rt}\leq
\varphi''(\uep_{r})\left(\left[\uep_{t}\right]^{2}+
\left[\uep_{rt}\right]^{2}\right)\leq
\varphi''(\uep_{r})\left[\uep_{rt}\right]^{2}+
c_{1},$$
hence
\begin{equation}
    I_{2}\leq\int_{1}^{\beta(t)}\varphi''(\uep_{r})
    \left[\uep_{rt}\right]^{2}\,dr+
    c_{2}.
    \label{est:int-i2}
\end{equation}

From (\ref{est:uept-glob}) we have also that
\begin{equation}
    I_{4}\leq c_{3}\beta'(t).
    \label{est:int-i4}
\end{equation}

Taking the time derivative of the Neumann boundary condition
(\ref{th:Q-nr}) we obtain that
$$0=\frac{\mbox{d}}{\mbox{d}t}\left[
\uep_{r}(\beta(t),t)\right]=\beta'(t)\uep_{rr}(\beta(t),t)+
\uep_{rt}(\beta(t),t).$$

From (\ref{est:ueprr-glob}) we have therefore that
$$|\uep_{rt}(\beta(t),t)|=\beta'(t)|\uep_{rr}(\beta(t),t)|
\leq c_{4}\beta'(t),$$
hence by  (\ref{est:uept-glob})
\begin{equation}
    I_{3}\leq c_{5}\beta'(t).
    \label{est:int-i3}
\end{equation}

Plugging (\ref{est:int-i2}), (\ref{est:int-i4}), and
(\ref{est:int-i3}) in (\ref{est:int-i}) we obtain that
$$\frac{\mbox{d}}{\mbox{d}t}\left(
\int_{1}^{\beta(t)}\left[\uep_{t}\right]^{2} dr\right)\leq
-\int_{1}^{\beta(t)}\varphi''(\uep_{r})
\left[\uep_{rt}\right]^{2}\,dr+c_{2}+c_{6}\beta'(t),$$
and therefore
$$\int_{0}^{t_{0}}\int_{1}^{\beta(t)}\varphi''(\uep_{r})
    \left[\uep_{rt}\right]^{2}dr\,dt\leq
    \int_{1}^{2}\left[
    \uep_{t}(r,0)\right]^{2}dr+c_{2}t_{0}+c_{6}(\beta(t_{0})-\beta(0))
	\leq c_{7}.$$

This proves (\ref{est:ueprt-int}).

In order to prove (\ref{est:ueprt}) through (\ref{est:ueprrt-int}), we
choose a cut-off function $\rho\in C^{\infty}(\re^{2})$ such that
$0\leq\rho(r,t)\leq 1$ for every $(x,t)\in\re^{2}$, and
\begin{eqnarray*}
	\rho(r,t) & = & \makebox[3em][l]{1}\mbox{in }\{(r,t)\in\qone: 1\leq
	r\leq\beta(t)-\delta\},  \\
	\rho(r,t) & = & \makebox[3em][l]{0}\mbox{in }\{(r,t)\in\qone:
	\beta(t)-\delta/2\leq
	r\leq\beta(t)\}.
\end{eqnarray*}

Then we set 
$$E(t):=\int_{1}^{\beta(t)}\rho^{2}[\uep_{rt}]^{2}\,dr,
\quad\quad\quad
F(t):=\int_{1}^{\beta(t)}\rho^{2}
\varphi''(\uep_{r})[\uep_{rrt}]^{2}\,dr,$$
$$G(t):=\int_{1}^{\beta(t)} \varphi''(\uep_{r})[\uep_{rt}]^{2}\,dr.$$

Taking the time derivative of $E(t)$ and integrating by parts we have
that
\begin{eqnarray*}
    E'(t) & = & 2\int_{1}^{\beta(t)}\rho\rho_{t}[\uep_{rt}]^{2}\,dr
	+2\int_{1}^{\beta(t)}\rho^{2}\uep_{rt}\uep_{rtt}\,dr   \\
     & = & 2\int_{1}^{\beta(t)}\rho\rho_{t}[\uep_{rt}]^{2}\,dr
	-2\int_{1}^{\beta(t)}(\rho^{2}\uep_{rt})_{r}\uep_{tt}\,dr,
\end{eqnarray*}
where we neglected the boundary terms because $u_{rt}=0$ when $r=1$,
and $\rho=0$ when $r=\beta(t)$. Computing the derivatives in the
right-hand side we end up with
\begin{eqnarray}
    E'(t)  & = & 2\int_{1}^{\beta(t)}\rho\rho_{t}[\uep_{rt}]^{2}\,dr
	\nonumber  \\
	 &  &  -4\int_{1}^{\beta(t)}\rho\rho_{r}	\varphi'''(\uep_{r})
	[\uep_{rt}]^{2}\uep_{rr}\,dr
	-4\int_{1}^{\beta(t)}\rho\rho_{r}\varphi''(\uep_{r})
	\uep_{rt}\uep_{rrt}\,dr
	\nonumber  \\
     &  &
	 -4\int_{1}^{\beta(t)}\rho\rho_{r}\frac{\varphi''(\uep_{r})}{r}
	 [\uep_{rt}]^{2}
	\,dr
	-2\int_{1}^{\beta(t)}\rho^{2}\varphi'''(\uep_{r})
	 \uep_{rt}\uep_{rr}\uep_{rrt}\,dr
	\nonumber  \\
     &  &
	 -2\int_{1}^{\beta(t)}\rho^{2}\varphi''(\uep_{r})
	 [\uep_{rrt}]^{2}\,dr
	 -2\int_{1}^{\beta(t)}\rho^{2}\frac{\varphi''(\uep_{r})}{r}
	 \uep_{rt}\uep_{rrt}\,dr
    \nonumber  \\
     & =: & I_{1}(t)+I_{2}(t)+I_{3}(t)+I_{4}(t)+I_{5}(t)+I_{6}(t)+I_{7}(t).
    \nonumber
\end{eqnarray}

Let us estimate separately the seven terms.  From now on all constants
depend on $\delta$, but are still independent on $\ep$.  

First of all we have that $I_{4}(t)\leq c_{8}G(t)$ and
$I_{6}(t)=-2F(t)$.

Due to the strict parabolicity in the interior and
(\ref{est:ueprr-glob}) we have that
$$I_{1}(t) \leq\int_{1}^{\beta(t)-\delta/2}|\rho_{t}|\cdot
[\uep_{rt}]^{2}\,dr\leq
c_{9}\int_{1}^{\beta(t)}\varphi''(\uep_{r})[\uep_{rt}]^{2}\,dr=
c_{9}G(t),$$
$$I_{2}(t) \leq c_{10}\int_{1}^{\beta(t)-\delta/2}
|\varphi'''(\uep_{r})|\cdot|\uep_{rr}|\cdot[\uep_{rt}]^{2}\,dr \leq
c_{11}G(t).$$

Moreover from inequality $2ab\leq\nu a^{2}+\nu^{-1}b^{2}$ we deduce that
$$I_{3}(t)\leq c_{12}\int_{1}^{\beta(t)}
\sqrt{\rho^{2} \varphi''(\uep_{r})[\uep_{rrt}]^{2}}\cdot
\sqrt{\varphi''(\uep_{r})[\uep_{rt}]^{2}}\,dr\leq
\frac{1}{3}F(t)+c_{13}G(t),$$
and similarly
$$I_{7}(t)\leq\frac{1}{3}F(t)+c_{14}G(t).$$

In the same way, exploiting once again the strict parabolicity in the
interior, and estimate (\ref{est:ueprr-glob}), we have that
$$I_{5}\leq c_{15}\int_{1}^{\beta(t)}\rho^{2}\varphi''(\uep_{r}) |
\uep_{rt}|\cdot|\uep_{rrt}|\,dr\leq \frac{1}{3}F(t)+c_{16}G(t).$$

Putting all together we obtain that
$E'(t)\leq -F(t)+c_{17}G(t)$,
hence
$$E(t)+\int_{0}^{t}F(s)\,ds\leq
E(0)+c_{17}\int_{0}^{t}G(s)\,ds.$$

The term $E(0)$ depends on the initial condition only.  The integral
of $G(t)$ is uniformly bounded because of (\ref{est:ueprt-int}).  It
follows that the left-hand side is bounded independently on $\ep$.
Thanks to the properties of $\rho(r,t)$, this proves (\ref{est:ueprt})
and (\ref{est:ueprrt-int}).  Finally, from (\ref{eq:v}) we have that
$$\varphi''(\uep_{r})|\uep_{rrr}|\leq|\uep_{rt}|+
|\varphi'''(\uep_{r})|\cdot|\uep_{rr}|^{2}+
\frac{\varphi''(\uep_{r})}{r}|\uep_{rr}|+
\frac{\varphi'(\uep_{r})}{r^{2}}.$$

From (\ref{est:ueprr-glob}) and the uniform parabolicity in the
interior it follows that
$$\int_{1}^{\beta(t)-\delta}|\uep_{rrr}|^{2}\,dr\leq c_{18}+
c_{19}\int_{1}^{\beta(t)-\delta}|\uep_{rt}|^{2}\,dr,$$
which proves (\ref{est:ueprrr-int}).

\subsection{Proof of Theorem~\ref{thm:Tep}}

Let us briefly sketch the outline of the proof, which is quite similar
to the proof of Theorem~\ref{thm:Qep}.  First of all we show that the
initial boundary value problem has a unique solution $\uep\in
C^{2,1}(\Tep)\cap C^{\infty}(\Tep\cap\intp(\T))$.  Then we set
$v=:\uep_{r}$, and $w:=\uep_{rr}$.  Both $v$ and $w$ belong to
$C^{0}(\Tep)\cap C^{\infty}(\Tep\cap\intp(\T))$.  It is easy to see
that $v$ is the solution of the following equation
\begin{equation}
	v_{t}=-\varphi''(v)v_{rr}-\varphi'''(v)v_{r}^{2}-
	\frac{\varphi''(v)}{r}v_{r}+\frac{\varphi'(v)}{r^{2}}
	\quad\quad\quad
	\forall (r,t)\in\intp(\Tep),
	\label{eq:T-v}
\end{equation}
with Dirichlet boundary conditions
\begin{equation}
	v(\beta(t_{0}-t),t)=v(\gamma(t_{0}-t),t)=1+\ep
	\quad\quad\quad
	\forall t\in[\ep,t_{0}],
	\label{eq:T-v-bcl}
\end{equation}
and initial datum 
\begin{equation}
	v(r,\ep)=1+\ep
	\quad\quad\quad
	\forall r\in[\beta(t_{0}-\ep),\gamma(t_{0}-\ep)].
	\label{eq:T-v-cauchy}
\end{equation}

In the same way $w$ is a solution in $\intp(\Tep)$ of equation (once
again the terms in $v$ are thought as coefficients)
\begin{eqnarray}
	w_{t} & = &
	-\varphi''(v)w_{rr}-3\varphi'''(v)w_{r}w-\varphi^{IV}(v)w^{3}
	\nonumber  \\
	 & & - \frac{\varphi'''(v)}{r}w^{2}-\frac{\varphi''(v)}{r}w_{r}
	 +2\frac{\varphi''(v)}{r^{2}}w-2\frac{\varphi'(v)}{r^{3}},
	\label{eq:T-w}
\end{eqnarray}
with Dirichlet boundary conditions
\begin{equation}
	w(\beta(t_{0}-t),t)=\uep_{rr}(\beta(t_{0}-t),t)
	\quad\quad\quad
	\forall t\in[\ep,t_{0}],
	\label{eq:T-w-bc-l}
\end{equation}
\begin{equation}
	w(\gamma(t_{0}-t),t)=\uep_{rr}(\gamma(t_{0}-t),t)
	\quad\quad\quad
	\forall t\in[\ep,t_{0}],
	\label{eq:T-w-bc-r}
\end{equation}
and initial datum
\begin{equation}
	w(r,\ep)=0
	\quad\quad\quad
	\forall r\in[\beta(t_{0}-\ep),\gamma(t_{0}-\ep)].
	\label{eq:T-w-cauchy}
\end{equation}

Then we show that $v$ and $w$ satisfy four sets of inequalities in
$\Tep$. 
\begin{itemize}
	\item The first pair of inequalities is
	\begin{equation}
		1+\ep\leq v(r,t)\leq 2+\varphi'(1)t.
		\label{est:T-ur-mp}
	\end{equation}
	
	From the second inequality in (\ref{hp:t01}) we have in particular
	that $t_{0}\leq 1/\varphi'(1)$.  Therefore (\ref{est:T-ur-mp})
	implies (\ref{est:Tep-ur}).

	\item The second inequality is that
	\begin{equation}
		v(r,t)\geq 1+\eta\left(\frac{t}{t_{0}}-(r-3)^{2}\right)
		\label{est:Tep-up}
	\end{equation}
	for a suitable constant $\eta>0$.  The term after $\eta$ is
	positive in $\intp(\T)$.  This implies (\ref{est:Tep-unif-farab}),
	hence also (\ref{est:Tep-unif-parab}).

	\item The third pair of inequalities is
	\begin{equation}
		v(r,t)\geq
		1+\ep+\left(b(t_{0}-t)-\sqrt{\ep}\right)(r-\beta(t_{0}-t))-
		\gamma_{0}\left(r-\beta(t_{0}-t)\right)^{2},
		\label{est:T-ueprr-ll}
	\end{equation}
	\begin{equation}
		v(r,t)\leq
		1+\ep+\left(b(t_{0}-t)+\sqrt{\ep}\right)(r-\beta(t_{0}-t))+
		\gamma_{0}\left(r-\beta(t_{0}-t)\right)^{2},
		\label{est:T-ueprr-lu}
	\end{equation}
	where $\gamma_{0}$ is the constant defined in (\ref{defn:gamma0}).
	Arguing as in the proof of Theorem~\ref{thm:Qep}, these
	inequalities yield 	(\ref{est:Tep-ueprr-l}). 
	
	In an analogous way we have that
	\begin{equation}
		v(r,t)\geq
		1+\ep+\left(c(t_{0}-t)+\sqrt{\ep}\right)(r-\gamma(t_{0}-t))-
		\gamma_{0}\left(r-\gamma(t_{0}-t)\right)^{2},
		\label{est:T-ueprr-rl}
	\end{equation}
	\begin{equation}
		v(r,t)\leq
		1+\ep+\left(c(t_{0}-t)-\sqrt{\ep}\right)(r-\gamma(t_{0}-t))+
		\gamma_{0}\left(r-\gamma(t_{0}-t)\right)^{2},
		\label{est:T-ueprr-ru}
	\end{equation}
	which imply (\ref{est:Tep-ueprr-r}). 

	\item Thanks to (\ref{est:Tep-ueprr-l}) and (\ref{est:Tep-ueprr-r}) we
	have an estimate on the values of $w$ at the boundary.  This is
	the starting point to prove the fourth set of inequalities
	\begin{equation}
		b(t_{0}-t)-\sqrt{\ep}-\gamma_{1}(r-\beta(t_{0}-t))\leq w(r,t)\leq
		b(t_{0}-t)+\sqrt{\ep}+\gamma_{1}(r-\beta(t_{0}-t)),
		\label{est:T-ueprr}
	\end{equation}
	\begin{equation}
		c(t_{0}-t)-\sqrt{\ep}-\gamma_{1}(\gamma(t_{0}-t)-r)\leq w(r,t)\leq
		c(t_{0}-t)+\sqrt{\ep}+\gamma_{1}(\gamma(t_{0}-t)-r),
		\label{est:T-ueprr-bis}
	\end{equation}
	where $\gamma_{1}$ is the constant defined in (\ref{defn:gamma1}).
	These imply (\ref{est:Tep-ueprr-ln}) and
	(\ref{est:Tep-ueprr-rn}).
\end{itemize}

We are now ready to proceed with the details.  Many steps of the proof
(for example the integral estimates) are analogous to the
corresponding steps in the proof of Theorem~\ref{thm:Qep}.  In these
cases we skip them, focussing only on what is different.  We also skip
the proofs of (\ref{est:T-ueprr-rl}), (\ref{est:T-ueprr-ru}),
(\ref{est:T-ueprr-bis}), which are analogous to the proofs of
(\ref{est:T-ueprr-ll}), (\ref{est:T-ueprr-lu}), and
(\ref{est:T-ueprr}), respectively.

\paragraph{Existence and maximum principle for space derivatives}

Let us take a function $\varphi_{\ep}\in C^{\infty}(\re)$ which
coincides with $\varphi$ in the interval $[1+\ep,3]$, and such that
$\varphi_{\ep}''(\sigma)\leq-\nu_{\ep}<0$ for every $\sigma\in\re$.
Equation (\ref{eq:T-eq}), with $\varphi_{\ep}$ instead of $\varphi$,
is strictly forward parabolic. Therefore by well know arguments the initial
boundary value problem admits a unique solution $\uep\in
C^{2,1}(\Tep)\cap C^{\infty}(\Tep\setminus\parbt)$.

It is easy to see that $z(r,t):=1+\ep$ is a subsolution of problem
(\ref{eq:T-v}) through (\ref{eq:T-v-cauchy}), both with
$\varphi_{\ep}$ and with $\varphi$.  This proves the lower bound in
(\ref{est:Tep-ur}). 

We claim that $z(r,t):=2+\varphi'(1)t$ is a supersolution of 
(\ref{eq:T-v}) through (\ref{eq:T-v-cauchy}). Indeed on the three
sides of the parabolic boundary of $\Tep$ we have that
$$z(r,t)\geq 2\geq 1+\ep=v(r,t).$$

In the interior we have that $z_{r}(r,t)=z_{rr}(r,t)=0$, hence (since 
$1\leq z(r,t)\leq 3$)
$$z_{t}(r,t)=\varphi'(1)\geq\varphi'(z(r,t))\geq
\frac{\varphi'(z(r,t))}{r},$$
which is exactly the required differential inequality (both with
$\varphi_{\ep}$ and with $\varphi$).

This completes the proof of (\ref{est:T-ur-mp}), hence of
(\ref{est:Tep-ur}).

\paragraph{Uniform strict parabolicity in the interior}

Let us choose $\eta>0$ small enough so that
\begin{equation}
	\eta\leq t_{0}\;(<1),
	\hspace{3em}
	\eta\leq\left(\frac{1}{t_{0}}+20\gamma_{2}\right)^{-1}\varphi'(3).
	\label{defn:T-eta}
\end{equation}

Let $z(r,t)$ denote the right-hand side of (\ref{est:Tep-up}).  We
claim that, under these restrictions on $\eta$, $z$ is a subsolution
of (\ref{eq:T-v}) through (\ref{eq:T-v-cauchy}), which implies
(\ref{est:Tep-up}).

\subparagraph{\emph{\textmd{Parabolic boundary}}}

When $r=\beta(t_{0}-t)$ or $r=\gamma(t_{0}-t)$ we have that
$$z(\beta(t_{0}-t),t)=z(\gamma(t_{0}-t),t)=1\leq 1+\ep=
v(\beta(t_{0}-t),t)=v(\gamma(t_{0}-t),t).$$

Since $\eta\leq t_{0}$, when $t=\ep$ we have that
$$z(r,\ep)=1+\eta\left(\frac{\ep}{t_{0}}-(r-3)^{2}\right)\leq
1+\frac{\eta\ep}{t_{0}}\leq 1+\ep=v(r,\ep).$$

\subparagraph{\emph{\textmd{Differential inequality}}}

We have that
$$z_{r}(r,t)=-2(r-3)\eta,
\hspace{2em}
z_{rr}(r,t)=-2\eta,
\hspace{2em}
z_{t}(r,t)=\eta/t_{0}.$$

In order to show that $z$ is a subsolution we have to prove that
$$\frac{\eta}{t_{0}}\leq 2\eta\varphi''(z)
-4\eta^{2}\varphi'''(z)(r-3)^{2}+ 2\eta\frac{\varphi''(z)}{r}(r-3)+
\frac{\varphi'(z)}{r^{2}}.$$

It is easy to show that $1\leq z(r,t)\leq 3$ for every $(r,t)\in\Tep$.
Therefore from the properties of $\varphi'$ we have that
$\varphi'(z)/r^{2}\geq(1/16)\varphi'(3)$.  All the other terms are
uniformly small when $\eta$ is small.  This shows that the
differential inequality is satisfied when $\eta$ is small enough, for
example as soon as $\eta$ fulfils the last inequality in
(\ref{defn:T-eta}).

\paragraph{Estimate on second derivatives at the moving boundary --
Subsolution} 

Let us prove (\ref{est:T-ueprr-ll}).  To this end we denote its
right-hand side by $z(r,t)$, and we show that it is a subsolution of
(\ref{eq:T-v}) through (\ref{eq:T-v-cauchy}).

\subparagraph{\emph{\textmd{Boundary $r=\beta(t_{0}-t)$}}}

In this case we have that
$$z(\beta(t_{0}-t),t)=1+\ep=v(\beta(t_{0}-t),t).$$

\subparagraph{\emph{\textmd{Boundary $r=\gamma(t_{0}-t)$}}}

From the explicit expressions (\ref{defn:beta-gamma}) we have that
\begin{eqnarray}
	z(\gamma(t_{0}-t),t) & = & 1+\ep+\left(
	b(t_{0}-t)-\sqrt{\ep}\right)\cdot 2\frac{\sqrt{t}}{\sqrt{t_{0}}}-
	4\gamma_{0}\frac{t}{t_{0}}
	\nonumber  \\
	 & \leq & 1+\ep+2\frac{\sqrt{t}}{\sqrt{t_{0}}}\left(
	 b(t_{0}-t)-2\gamma_{0}\frac{\sqrt{t}}{\sqrt{t_{0}}}\right).
	\label{eq:z-gamma}
\end{eqnarray}

Exploiting estimate (\ref{est:b-sqrt}), and the fact that
$\gamma_{0}\geq 1/2$ and $t_{0}\leq 1$, we have that
\begin{equation}
	b(t_{0}-t)\leq \sqrt{t}\leq
	2\gamma_{0}\sqrt{t}\leq 2\gamma_{0} \frac{\sqrt{t}}{\sqrt{t_{0}}}.
	\label{eq:stima}
\end{equation}

From (\ref{eq:z-gamma}) and (\ref{eq:stima}) we conclude that
$z(\gamma(t_{0}-t),t)\leq 1+\ep=v(\gamma(t_{0}-t),t)$.

\subparagraph{\emph{\textmd{Boundary $t=\ep$}}}

From (\ref{est:b-sqrt}) we have that $b(t_{0}-\ep)\leq\sqrt{\ep}$,
hence
\begin{eqnarray*}
	z(r,\ep) & = &
	1+\ep+\left(b(t_{0}-\ep)-\sqrt{\ep}\right)(r-\beta(t_{0}-\ep))-
	\gamma_{0}\left(r-\beta(t_{0}-\ep)\right)^{2} \\
	 & \leq & 1+\ep\ =\ v(r,\ep).
\end{eqnarray*}

Note that the term with $\sqrt{\ep}$ instead of $\ep$ in the
right-hand side of (\ref{est:T-ueprr-ll}) is essential in this point
of the proof.

\subparagraph{\emph{\textmd{Differential inequality}}}

Let us set for simplicity $x:=r-\beta(t_{0}-t)$. Then 
we have that
$$z_{r}(r,t)=b(t_{0}-t)-\sqrt{\ep}-2\gamma_{0}x,
\hspace{3em}
z_{rr}(r,t)=-2\gamma_{0},$$
$$z_{t}(r,t)=b(t_{0}-t)\beta'(t_{0}-t)-\beta'(t_{0}-t)\sqrt{\ep}-
(b'(t_{0}-t)+2\gamma_{0}\beta'(t_{0}-t))x.$$

Plugging these expressions in (\ref{eq:T-v}), and exploiting the
equation defining $b(t_{0}-t)$ as we did in (\ref{eq-beta-appl}), we
end up with the following differential inequality
\begin{eqnarray*}
	\lefteqn{\hspace{-5em} \beta'(t_{0}-t)\sqrt{\ep}
	+(b'(t_{0}-t)+2\gamma_{0}\beta'(t_{0}-t))x \ \geq} \\
	 & \geq & -2\gamma_{0}\varphi''(z)  \\
	 \mbox{\hspace{5em}} &  & +\varphi'''(z)\left(b(t_{0}-t)-
	 \sqrt{\ep} -2\gamma_{0}x\right)^{2}
	 -\varphi'''(1)b^{2}(t_{0}-t)    \\
	 &  & +\frac{\varphi''(z)}{r}
	 \left(b(t_{0}-t)-\sqrt{\ep}-2\gamma_{0}x\right)  \\
	 &  & +\frac{\varphi'(1)}{\beta^{2}(t_{0}-t)}-
	 \frac{\varphi'(z)}{r^{2}}.  
\end{eqnarray*}

This inequality has to be proved for all $(r,t)\in\Tep$ such that
$1\leq z(r,t)\leq 3$.  This can be done arguing as we did in the
corresponding step of the proof of Theorem~\ref{thm:Qep} (and using
from time to time that $\sqrt{\ep}\geq\ep$).

\paragraph{Estimate on second derivatives at the moving boundary --
Supersolution} 

In order to prove (\ref{est:T-ueprr-lu}) it is enough to show that its
right-hand side is a supersolution of (\ref{eq:T-v}) through
(\ref{eq:T-v-cauchy}).  In this case the inequalities on the three
pieces of the parabolic boundary of $\Tep$ are trivial.  The proof of
the differential inequality is analogous to the case of the
supersolution. We skip the details.

\paragraph{Global estimate on second derivatives -- Subsolution} 

Let $z(r,t)$ denote the left-hand side of (\ref{est:T-ueprr}).  We
claim that it is a subsolution of (\ref{eq:T-w}) through
(\ref{eq:T-w-cauchy}).

\subparagraph{\emph{\textmd{Boundary $r=\beta(t_{0}-t)$}}}

From (\ref{est:Tep-ueprr-l}) we have that
$$z(\beta(t_{0}-t),t)=b(t_{0}-t)-\sqrt{\ep}\leq
\uep_{rr}(\beta(t_{0}-t),t)=w(\beta(t_{0}-t),t).$$

\subparagraph{\emph{\textmd{Boundary $r=\gamma(t_{0}-t)$}}}

Exploiting estimates (\ref{est:b-sqrt}), (\ref{est:c-sqrt}), and
(\ref{est:Tep-ueprr-r}), and the fact that $\gamma_{1}\geq 1$ and
$t_{0}\leq 1$, we have that
$$z(\gamma(t_{0}-t),t)=b(t_{0}-t)-\sqrt{\ep}-
2\gamma_{1}\frac{\sqrt{t}}{\sqrt{t_{0}}}\leq\sqrt{t}-\sqrt{\ep}-
2\gamma_{1}\frac{\sqrt{t}}{\sqrt{t_{0}}}\leq
-\sqrt{t}-\sqrt{\ep}\leq$$
$$\leq c(t_{0}-t)-\sqrt{\ep}\leq
\uep_{rr}(\gamma(t_{0}-t),t)=w(\gamma(t_{0}-t),t).$$

\subparagraph{\emph{\textmd{Boundary $t=\ep$}}}

From (\ref{est:b-sqrt}) we have that $b(t_{0}-\ep)\leq\sqrt{\ep}$,
hence
$$z(r,\ep)=b(t_{0}-\ep)-\sqrt{\ep}-\gamma_{1}(r-\beta(t_{0}-\ep)) 
\leq 0=w(r,\ep).$$

\subparagraph{\emph{\textmd{Differential inequality}}}

After computing the derivatives and changing the signs we have to
prove that 
$$b'(t_{0}-t)+\gamma_{1}\beta'(t_{0}-t)\geq
-3\gamma_{1}\varphi'''(v)z+
\varphi^{IV}(v)z^{3}+\frac{\varphi'''(v)}{r}z^{2}-
\gamma_{1}\frac{\varphi''(v)}{r}
-2\frac{\varphi''(v)}{r^{2}}z+2\frac{\varphi'(v)}{r^{3}}.$$

Note that this is exactly the same inequality satisfied by
supersolutions in $\qone$ (since we reversed the time, subsolutions
correspond to supersolutions and vice versa).  The proof is of course
completely analogous.

\paragraph{Global estimate on second derivatives -- Supersolution} 

Let $z(r,t)$ denote the right-hand side of (\ref{est:T-ueprr}).  We
claim that $z$ is a supersolution of (\ref{eq:T-w}) through
(\ref{eq:T-w-cauchy}).

\subparagraph{\emph{\textmd{Boundary $r=\beta(t_{0}-t)$}}}

From (\ref{est:Tep-ueprr-l}) we have that
$$z(\beta(t_{0}-t),t)=b(t_{0}-t)+\sqrt{\ep}\geq
\uep_{rr}(\beta(t_{0}-t),t)=w(\beta(t_{0}-t),t).$$

\subparagraph{\emph{\textmd{Boundary $r=\gamma(t_{0}-t)$}}}

We have that
$$z(\gamma(t_{0}-t),t)=b(t_{0}-t)+\sqrt{\ep}+
2\gamma_{1}\frac{\sqrt{t}}{\sqrt{t_{0}}}\geq\sqrt{\ep}.$$

Moreover, from (\ref{est:Tep-ueprr-r}) and the fact that
$c(t_{0}-t)\leq 0$, we have that
$$w(\gamma(t_{0}-t),t)=\uep_{rr}(\gamma(t_{0}-t),t)\leq
c(t_{0}-t)+\sqrt{\ep}\leq\sqrt{\ep}.$$

It follows that $z(\gamma(t_{0}-t),t)\geq w(\gamma(t_{0}-t),t)$.

\subparagraph{\emph{\textmd{Boundary $t=\ep$}}}

We have trivially that
$z(r,\ep)\geq 0=w(r,\ep)$.

\subparagraph{\emph{\textmd{Differential inequality}}}

We have to prove the same inequality satisfied by the
corresponding subsolution in $\qone$. The proof is the same.

\setcounter{equation}{0}
\section{Passing to the limit}\label{sec:limit}

The basic tool is the following compactness result.  We omit the
proof, for which we refer to \cite[Lemma~3.1]{GG-LocSol}.  The key
point in that (\ref{hp:L3}) through (\ref{hp:L1}) yield a uniform
bound on the norm of $\fep$ in the H\"{o}lder space
$C^{1/2,1/4}([r_{1},r_{2}]\times[t_{1},t_{2}])$.

\begin{lemma}\label{lemma:conv}
	Let $[r_{1},r_{2}]\times[t_{1},t_{2}]$ be a rectangle, and let
	$\ep_{0}>0$.  For every $\ep\in(0,\ep_{0})$, let $\fep\in
	C^{1}([r_{1},r_{2}]\times[t_{1},t_{2}])$.  Let us assume that
	there exists $M\in\re$ such that
    \begin{eqnarray}
        |\fep(r_{1},t_{1})| & \leq & M\hspace{2em}\forall
		\ep\in(0,\ep_{0}),
        \label{hp:L3}\\
        \int_{r_{1}}^{r_{2}}[\fep_{r}(r,t)]^{2}\,dr & \leq & 
		M\hspace{2em}\forall t\in[t_{1},t_{2}],\ \forall \ep\in(0,\ep_{0}),
        \label{hp:L2}  \\
		\int_{t_{1}}^{t_{2}}\int_{r_{1}}^{r_{2}}
		[\fep_{t}(r,t)]^{2}\,dr\,dt & \leq & M\hspace{2em}\forall
		\ep\in(0,\ep_{0}).
        \label{hp:L1}  
    \end{eqnarray}
    
    Then the family $\{\fep\}$ is relatively compact in
    $C^{0}([r_{1},r_{2}]\times[t_{1},t_{2}])$. 
	\qed
\end{lemma}

\subsection{Proof of Theorem~\ref{thm:q1}}

Let $t_{0}$ and $u_{0}$ be as in Theorem~\ref{thm:Qep}, and let
$\{\uep\}_{\ep\in(0,1)}$ be the family of solutions produced by that
theorem.  We claim that $\uep$ uniformly converges, as $\ep\to 0^{+}$,
to a limit $u$, which satisfies the initial condition
\begin{equation}
	u(r,0)=u_{0}(r)
	\quad\quad
	\forall r\in[1,2]
	\label{eq:u-cauchy}
\end{equation}
and all the conditions required in Theorem~\ref{thm:q1}.

\subparagraph{\emph{\textmd{Uniqueness}}}

Let us assume that a function $u\in C^{2,1}(\qone)$ satisfies equation
(\ref{eq:PM-rad}) for every $(r,t)\in\qone$, the Neumann boundary
condition (\ref{eq:nbc-l}) for every $t\in[0,t_{0}]$, the Neumann
boundary condition (\ref{g1:ur}), and the initial condition
(\ref{eq:u-cauchy}).  Let us assume also that $0\leq u_{r}(r,t)\leq 1$
for every $(r,t)\in\qone$.  Under this condition equation
(\ref{eq:PM-rad}) is degenerate parabolic, hence the solution $u$ is
unique.

This shows that the limit problem provides a unique characterization
of the possible limits.  Therefore in what follows we can limit
ourselves to show that $\uep$ converges to the solution of this
problem up to subsequences.

\subparagraph{\emph{\textmd{Convergence of $\uep$}}}

Estimates (\ref{est:Q-ur}) and (\ref{est:uept-glob}) provide a uniform
bound on the Lipschitz constant of $\uep$.  Moreover the functions
$\uep$ are equi-bounded due to the initial condition.  From the
classical Ascoli's Theorem it follows that (up to subsequences, which we
don't relabel) $\uep$ uniformly converges in $\qone$ to a continuous
function $u$, satisfying of course the initial condition
(\ref{eq:u-cauchy}).

\subparagraph{\emph{\textmd{Convergence of $\uep_{r}$}}}

Let us consider any rectangle
$[r_{1},r_{2}]\times[t_{1},t_{2}]\subseteq\qone\setminus\Gamma_{1}$, and
let us apply Lemma~\ref{lemma:conv} in this rectangle with
$\fep=\uep_{r}$.  In the rectangle assumption (\ref{hp:L3}) follows
follows from (\ref{est:Q-ur}), assumption (\ref{hp:L2}) follows from
(\ref{est:ueprr-glob}), and assumption (\ref{hp:L1}) follows from
(\ref{est:ueprt-int}) and the strict parabolicity in the interior.  We
obtain that 
$$\uep_{r}\to u_{r} \quad\quad
\mbox{uniformly in }[r_{1},r_{2}]\times[t_{1},t_{2}].$$

Taking the union over all such rectangles, we easily conclude that
$$\uep_{r}(r,t)\to u_{r}(r,t)
\quad\quad
\forall(r,t)\in\qone\setminus\Gamma_{1}.$$

We claim that $u_{r}$ can be continuously extended to the whole
$\qone$ by setting $u_{r}(r,t)=1$ for every $(r,t)\in\Gamma_{1}$.  In
order to prove the continuity of this extension, let us consider
$(r_{1},t_{1})=(\beta(t_{1}),t_{1})\in\Gamma_{1}$ and
$(r_{2},t_{2})\in\qone\setminus\Gamma_{1}$.  From
(\ref{est:ueprr-glob}) we have that
$$|\uep_{r}(r,t)-(1+\ep)|=|\uep_{r}(r,t)-\uep_{r}(\beta(t),t)|\leq
M_{4}|r-\beta(t)|.$$

Letting $\ep\to 0^{+}$ we obtain that
$$|u_{r}(r,t)-1|\leq M_{4}|r-\beta(t)|
\quad\quad\quad
\forall (r,t)\in\qone\setminus\Gamma_{1}.$$

Therefore we have that
\begin{eqnarray*}
    |u_{r}(r_{2},t_{2})-u_{r}(r_{1},t_{1})| & = &
    |u_{r}(r_{2},t_{2})-1|  \\
     & \leq & M_{4}|r_{2}-\beta(t_{2})|  \\
     & \leq & M_{4}\left(|r_{2}-r_{1}|+
     |\beta(t_{1})-\beta(t_{2})|\right).
\end{eqnarray*}

Due to the continuity of $\beta(t)$, the right-hand side is small when
$(r_{2},t_{2})$ is close to $(r_{1},t_{1})$.  This proves the
continuity of $u_{r}$ up to $\Gamma_{1}$.  Of course $u_{r}$ satisfies
the Neumann boundary conditions in $r=1$ and in $r=\beta(t)$.

\subparagraph{\emph{\textmd{Convergence of $\uep_{rr}$}}}

Let us consider any rectangle
$[r_{1},r_{2}]\times[t_{1},t_{2}]\subseteq\qone\setminus\Gamma_{1}$,
and let us apply Lemma~\ref{lemma:conv} in this rectangle with
$\fep=\uep_{rr}$.  In the rectangle assumption (\ref{hp:L3}) follows
follows from (\ref{est:ueprr-glob}), assumption (\ref{hp:L2}) follows
from (\ref{est:ueprrr-int}), and assumption (\ref{hp:L1}) follows from
(\ref{est:ueprrt-int}).  We obtain that $\uep_{rr}\to u_{rr}$
uniformly in $[r_{1},r_{2}]\times[t_{1},t_{2}]$.

Taking the union over all such rectangles, we easily conclude that
$$\uep_{rr}(r,t)\to u_{rr}(r,t)
\quad\quad
\forall(r,t)\in\qone\setminus\Gamma_{1}.$$

Let us extend $u_{rr}$ to the whole $\qone$ by setting
$u_{rr}(r,t)=b(t)$ for every $(r,t)\in\Gamma_{1}$.  In order to prove
the continuity of this extension, let us consider
$(r_{1},t_{1})=(\beta(t_{1}),t_{1})\in\Gamma_{1}$ and
$(r_{2},t_{2})\in\qone\setminus\Gamma_{1}$.  Passing
(\ref{est:Q-ueprr-n}) to the limit we obtain that
$$|u_{rr}(r,t)-b(t)|\leq M_{3}|r-\beta(t)| \quad\quad\quad
\forall (r,t)\in\qone\setminus\Gamma_{1},$$
hence
\begin{eqnarray*}
    |u_{rr}(r_{2},t_{2})-u_{rr}(r_{1},t_{1})| & = &
    |u_{rr}(r_{2},t_{2})-b(t_{1})|  \\
     & \leq & |u_{rr}(r_{2},t_{2})-b(t_{2})|+|b(t_{2})-b(t_{1})|  \\
     & \leq & M_{3}|r_{2}-\beta(t_{2})|+|b(t_{2})-b(t_{1})|  \\
     & \leq & M_{3}\left(|r_{2}-r_{1}|+
	 |\beta(t_{2})-\beta(t_{1})|\right)+
     |b(t_{2})-b(t_{1})|.
\end{eqnarray*}

The conclusion easily follows from the continuity of $\beta(t)$ and
$b(t)$.

\subparagraph{\emph{\textmd{Convergence of $\uep_{t}$}}}

Since $\uep_{t}$ is related to $\uep_{r}$ and $\uep_{rr}$ by
(\ref{eq:PM-rad}), the convergence of $\uep_{t}$ to a continuous
function in $\qone$ follows from the convergence of $\uep_{r}$ and
$\uep_{rr}$ and the continuity of their limits.  The limit of
$\uep_{t}$ is of course $u_{t}$.  This completes the proof that $u$ is
of class $C^{2,1}$ and solves equation (\ref{eq:PM-rad}).

\subsection{Proof of Theorem~\ref{thm:q2}}

We have already seen in section~\ref{sec:back} that
Theorem~\ref{thm:q2} is equivalent to Theorem~\ref{thm:T}. In turn,
the proof of Theorem~\ref{thm:T} follows by passing to the limit the
solutions $\uep$ provided by Theorem~\ref{thm:Tep}.

The argument is quite similar to the proof of Theorem~\ref{thm:Qep}.
Thanks to some of the estimates of Theorem~\ref{thm:Tep} we can apply
Lemma~\ref{lemma:conv} and pass to the limit in any rectangle
$[r_{1},r_{2}]\times[t_{1},t_{2}]\subseteq\intp(\T)$, where $\uep$ is
defined for every $\ep$ small enough.  Thanks to the remaining
estimates of Theorem~\ref{thm:Tep} we can show that the limit $u$ and
its derivatives $u_{r}$, $u_{rr}$, $u_{t}$ can be continuously
extended up to the boundary.  We skip the details.

\label{NumeroPagine}

\end{document}